\documentclass[11pt]{amsart}
\usepackage{amsfonts}
\usepackage{cases}
\usepackage{mathrsfs}
\usepackage{bbm}
\usepackage{amssymb}
\usepackage{txfonts}
\usepackage{amscd}
\usepackage{dsfont}
\usepackage{amsfonts,latexsym,amsmath,amsthm,amsxtra,mathdots}
\usepackage[bookmarks,colorlinks]{hyperref}
\usepackage[all,cmtip]{xy}
\RequirePackage{amsmath} \RequirePackage{amssymb}
\usepackage{color}
\usepackage{colordvi}
\usepackage{multicol}
\usepackage{tikz}
\usepackage{hyperref}
\usepackage{graphicx}
\usepackage{amsmath}
\usepackage{amsmath,amscd}
\usetikzlibrary{matrix,arrows}

\def\-{^{-1}}

\newcommand{\delete}[1]{}

    \theoremstyle{plain}

\theoremstyle{plain}
\newtheorem{thm}{Theorem}[section] 

\theoremstyle{definition}
\newtheorem{defn}[thm]{Definition} 
\newtheorem{lem}[thm]{Lemma}
\newtheorem{prop}[thm]{Proposition}
\newtheorem{cor}[thm]{Corollary}
\newtheorem{rem}[thm]{Remark}

    \numberwithin{equation}{section}

\def\Proof{\noindent{\bf Proof}\quad}
\def\qed{\hfill$\square$\smallskip}

\begin{document}

\title{The Farrell-Jones Conjecture for the solvable Baumslag-Solitar groups}
\author{Tom Farrell and Xiaolei Wu}

\begin{abstract}
In this paper, we
prove the Farrell-Jones Conjecture for the solvable Baumslag-Solitar groups with coefficients in an additive category. We also extend our results to groups of the form, Z[1/p] semidirect product with  any virtually cyclic group, where p is a prime number.

\end{abstract}
\footnotetext{Date: November, 2012}



\maketitle


\section{Introduction}

In the Farrell and  Linnell paper \cite{FL}, they proved that if the fibered isomorphism conjecture is true for all nearly crystallographic groups, then it is true for all virtually solvable groups. However, they were not able to verify the fibered isomorphism conjecture for all nearly crystallographic groups. In particular, they pointed out that the fibered isomorphism conjecture has not been verified for the group ${\mathbb{Z} {[{\frac{1}{2}} ]}} \rtimes_\alpha \mathbb{Z}$, where $\alpha$ is multiplication by $2$. Note this group is isomorphic to the Baumslag-Solitar group $BS(1,2)$. Recall that the Baumslag-Solitar group $ BS(m,n)$ is defined by $\langle a,b~|~ba^mb^{-1} = a^n \rangle$ and all the solvable ones are isomorphic to $BS(1,d)$. Note that $BS(m,n) \cong BS(n,m) \cong BS(-m,-n)$. Using new technology developed by Bartels, L\"uck and Reich in \cite{BL1},\cite{BL2}, \cite{BLR}, \cite{BLR2}, we prove the following result.
\begin{thm} \label{mth}
The K-theoretic and L-theoretic Farrell-Jones Conjecture is true for all solvable Baumslag-Solitar groups with coefficients in an additive category.

\end{thm}

Note the truth of the Farrell-Jones conjecture with coefficients in an additive category implies the fibered isomorphism conjecture.  For more information about the Farrell-Jones conjecture and its fibered version, see for example \cite{FJ}. For the precise formulation and discussion of the Farrell-Jones conjecture
with coefficients in an additive category, see for example \cite{BFL}, \cite{BLR2}.  For our convenience, we will first prove the Farrell-Jones conjecture for $BS(1,d)$, where $d >1$. The same proof applies to the $d < -1$  case. Note that $BS(1,1)$ and $B(1,-1)$ are the fundamental groups of the torus and Klein bottle respectively. And the Farrell-Jones conjecture is known for those groups. The authors want to point out here that our current method can not be applied to all Baumslag-Solitar groups. For example, we do not know whether the Farrell-Jones conjecture is true for the group $BS(2,3)$. For convenience we will denote the group $BS(1,d)$ as ${\mathbb{Z} {[{\frac{1}{d}} ]}} \rtimes_\alpha \mathbb{Z}$ for the rest of the paper.
\begin{rem}
C. Wegner generalized the method in this paper and proved the Farrell-Jones conjecture for all virtually solvable groups in \cite{W2}. Using Wegner's result, we proved the Farrell-Jones conjecture for all Baumslag-Solitar groups in \cite{FW3}. Independently, G. Gandini, S. Meinert and H. R\"uping proved the Farrell-Jones conjecture for the fundamental group of any graph of abelian groups in \cite{GMR}, which includes all Baumslag-Solitar groups.
\end{rem}

Our strategy is to show that ${\mathbb{Z} {[{\frac{1}{d}} ]}} \rtimes_\alpha \mathbb{Z}$ is in fact a Farrell-Hsiang group, defined by  Bartels and L\"uck in \cite{BL1}. The main difficulty is to find a suitable flow space, ours  is a horizontal subspace of Bartels and L\"uck's flow space in \cite{BL2}. The horizontal flow space allows us to exploit some negative curvature present in the solvable Baumslag-Solitar groups. Once our flow space is defined, we prove in Lemma \ref{l3} that it has a lot of good properties, which implies that it has a
$\mathcal{VC}yc$-cover by results of Bartels, L\"uck and Reich in \cite{BLR}. It should be mentioned that our results rely heavily on the existence of such a cover. In the proof of our main theorem, we also use a Proposition of  Bartels, L\"uck and Reich from \cite{BLR2}, Proposition 5.3. In the last section, we extend our results to any group of the form $\mathbb{Z} {[{\frac{1}{p}} ]} \rtimes C$, where $p$ is any prime number, C is any virtually cyclic group.

In this paper, without further assumption, when we write ${\mathbb{Z} {[{\frac{1}{d}} ]}} \rtimes_\alpha \mathbb{Z}$, $\alpha$ is always multiplication by $d$, and we will assume $d$ is a fixed integer bigger than $1$. When we write Farrell-Jones conjecture, we mean the Farrell-Jones conjecture with coefficients in an additive category. Let G be a (discrete) group acting on a space X. We say the action is proper if for any $x \in X$ there is an open neighborhood U of $x$ such that $\{g \in G ~|~gU \bigcap U \neq {\emptyset}\}$ is finite.

\textbf{Acknowledgements.} This research was in part supported by the National Science Foundation. We would like to thank Arthur Bartels, Robert Bieri, Zhi Qi and Adrian Vasiu for helpful discussions. The second author also want to take this opportunity to thank his advisor, Tom Farrell, for introducing him to this wonderful field. We also thank the referees for many suggestions on how to improve the
readability of this paper and pointing out many typos in the previous version of the paper.

\section{A Model for $E({\mathbb{Z} {[{\frac{1}{d}} ]}} \rtimes_\alpha \mathbb{Z})$} \label{model}

In this section, we give a model for $E({\mathbb{Z} {[{\frac{1}{d}} ]}} \rtimes_\alpha \mathbb{Z})$, a contractible space
with free, proper and discontinuous $\Gamma$ action, where $\Gamma$ is the group ${\mathbb{Z} {[{\frac{1}{d}} ]}} \rtimes_\alpha \mathbb{Z}$. We also put a metric on it, such that $\Gamma$ acts isometrically on
$E({\mathbb{Z} {[{\frac{1}{d}} ]}} \rtimes_\alpha \mathbb{Z})$. Most part of the material in this section has been well studied before, see for example  \cite{ECH}, section 7.4, and \cite{FM}.

Let X be  $S^1 \times [0,1] / (z,0) \thicksim (z^d, 1)$. Its fundamental group is $\langle x,k  ~|~ kxk^{-1} ~=~ x^d\rangle$, which is isomorphic to our group  ${\mathbb{Z} {[{\frac{1}{d}} ]}} \rtimes_\alpha \mathbb{Z}$. Since its universal cover is contractible, it is a model for $E({\mathbb{Z} {[{\frac{1}{d}} ]}} \rtimes_\alpha \mathbb{Z})$.

\begin{figure}[h]
		\begin{center}
		\includegraphics[width=1.0\textwidth]{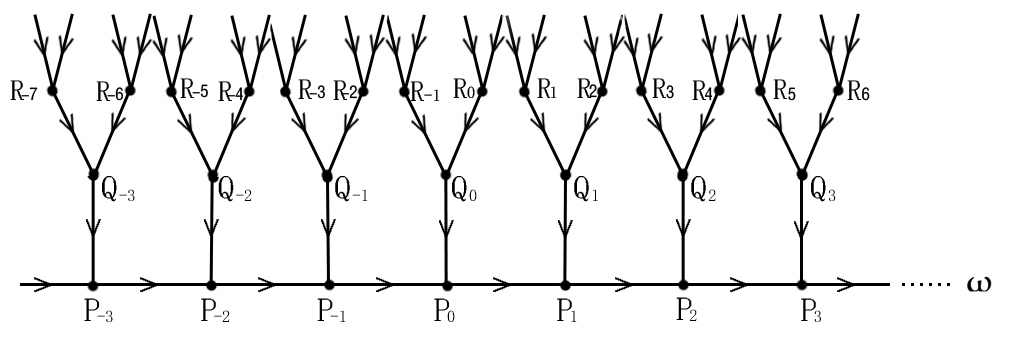}
		\end{center}
        \caption{$T_2$ : Oriented Infinite Regular 3-valent Tree}
        \label{tree}
	\end{figure}

 Let $T_d$ be the oriented infinite, regular, $(d+1)$-valent tree with edge length $1$ (see Figure \ref{tree}). At each vertex there are $d$ incoming and $1$ outgoing edges. There is a natural clockwise (from left to right) order on the $d$ incoming edges when we embed $T_d$ in the upper half plane with a specified horizontal line as shown in Figure \ref{tree}. We will denote the specified horizontal line by $L_0$, $L_0$ is the line $\ldots P_{-2}P_{-1}P_{0}P_{1}P_{2}\ldots$ in figure \ref{tree}. It is not too hard to show that the universal cover of X is $T_d \times \mathbb{R}$. However, it is not easy to figure out how $\Gamma$ acts on $T_d \times \mathbb{R}$. So what we do here is give an action of $\Gamma$ on $T_d \times \mathbb{R}$, and show that ${(T_d \times \mathbb{R})}  /\Gamma$ is homeomorphic to X.

We first assume $d$ is a prime number. We can embed $\Gamma$ into $GL_2(\mathbb{Q})$ by mapping $(x,k) \in {\mathbb{Z} {[{\frac{1}{d}} ]}} \rtimes_\alpha \mathbb{Z}$  to $\left(
\begin{array}{cc}
    d^{-k} & x \\
   0 & 1 \\
    \end{array}\right)$. There is a natural action of $GL_2(\mathbb{Q})$ on this tree. We explain briefly how $\Gamma$ acts on it. For more information, see for example \cite{Se}, page 69 - 78. The action of $\Gamma$ on $T_d$ when $d$ is not a prime can  be induced from the case when $d$ is a prime, details will be put in the Appendix.

Note that $\Gamma$ is generated by  $\left(
\begin{array}{cc}
    d^{-1} & 0 \\
   0 & 1 \\
    \end{array}\right)$ and $\left(
\begin{array}{cc}
    1 & d^{-m} \\
   0 & 1 \\
    \end{array}\right)$ , so we only need to show how these elements act on $T_d$. Choose  the vertex $P_0$ in $L_0$ as the basepoint. Denote the infinite point along $L_0$ towards the positive direction as $\omega$; cf. Figure \ref{tree}. For any $z \in T_d$, there is a unique geodesic starting from $z$ and moving towards $\omega$; denote it as $[x, \omega)$. Define the action  of $\left(
\begin{array}{cc}
    d^{-1} & 0 \\
   0 & 1 \\
    \end{array}\right)$  on the tree by translation along $L_0$, moving along the line $L_0$ towards the positive direction by $1$ units. For example in the case $d =2$ as in Figure \ref{tree}, $\left(
\begin{array}{cc}
    d^{-1} & 0 \\
   0 & 1 \\
    \end{array}\right)$  acts on the tree by translation along $L_0$, moving right along the line $L_0$ by $1$ units. Hence $P_n$ will be translated to $P_{n+1}$, and every vertex or edge growing on the tree with root $P_n$ will also be translated to the corresponding vertex or edge growing from root $P_{n+1}$.  We first give some notation and terminology before we  explain how $\beta_m = {\left(
\begin{array}{cc}
    1 & d^{-m}  \\
   0 & 1 \\
    \end{array}\right)}$ acts on $T_d$. We define a Busemann function $f_d:T_d\rightarrow\mathbb{R}$, by mapping the line $L_0$ isometrically to $\mathbb{R}$ with $P_0$ mapped to $0$ and $\omega$ to $-\infty$ (orientation reversed); for an arbitrary point $Q$, define $f_d(Q) = f_d(P) + d(P,Q)$, where $P$ is the closest point in the line $L_0$ from $Q$. For example in the case $d =2$ as in Figure \ref{tree}, we map $L_0$ to the real line $\mathbb{R}$ with $f_d(P_n) = -n$. Then for example $f_d(R_{1}) = 1$ since the closest point to $R_1$ is $P_1$ and $d(R_1,P_1) = 2$. Let $H_n = f_d^{-1}(n)$ which we call a horosphere in $T_d$ with center $\omega$. Likewise, let $B_n = f_d^{-1}((-\infty,n])$ be the corresponding horoball; cf. [\cite{BGS}, p23] for this terminology. And if $z \in H_n$, let $T_d(z)$ be the subtree in $T_d$ rooted at $z$ and growing outside $B_n$. Finally, for each $z \in H_n$ and $l \in {\mathbb{Z}}^+$, let
      $$ S(z,l) ~=~\{z' \in T_d(z) ~|~ d(z,z') =l\}.$$

\begin{rem}\label{bus}
These definitions of $f_d$, $H_n$, $B_n$ and $S(z,l)$ for the tree $T_d$ are valid for any positive integer $d$ (not just for primes).
\end{rem}

 With the terminology above, we describe the key features of the $\beta_m$ action on $T_d$  as follows (keep in mind $\beta_m^d = \beta_{m-1}$ and  $\left(
\begin{array}{cc}
    d^{-1} & 0 \\
   0 & 1 \\
    \end{array}\right)$ conjugates $\beta_m$ to $ \beta_{m+1}$):
 \begin{itemize}
\item[(i)] $\beta_m$ fixes each point in $B_{-m}$;
\item[(ii)] For each $z \in H_{-m}$ and $l \in {\mathbb{Z}}^+$, $\beta_m$ leaves $S(z,l)$ invariant and cyclically permutes its members. In particular, for each $z \in H_{-m}$, if we label elements in $S(z,1)$ by $(1,2,\cdots, d)$ with order preserved, then $\beta_m$ maps $(1,2,\cdots, d-1,d)$ to $(2,3,\cdots,d,1)$.
\end{itemize}

Note that $P_n$ has isotropy $\{ \left(
\begin{array}{cc}
    1 & d^{-n}s \\
   0 & 1 \\
    \end{array}\right)~|~ s \in \mathbb{Z}\}$. One of the good things about this action is that it fixes the infinite point $\omega$. We also define the action of $\Gamma$ on $\mathbb{R}$ by ${\left(
\begin{array}{cc}
    d^{-k} & x \\
   0 & 1 \\
    \end{array}\right)} ~w = {d^{-k}w + x}$, for $w \in \mathbb{R}$.

With these preparations, we define the action of $\Gamma$ on $T_d \times \mathbb{R}$ to be the diagonal action, i.e. $g  (z,w) = (gz,gw)$, for $g \in \Gamma, ~(z,w) \in T_d \times \mathbb{R}$. It is not hard to check that $\Gamma$ acts freely, properly and discontinuously on
$T_d \times \mathbb{R}$.  A fundamental domain is $P_0P_1\times [0,1]$, while very luckily ${(T_d \times \mathbb{R})}/\Gamma$ is homeomorphic to our previous space X. So $T_d \times \mathbb{R}$ is a model for $E({\mathbb{Z} {[{\frac{1}{d}} ]}} \rtimes_\alpha \mathbb{Z})$.

\begin{rem}{\label{act}}

Let $G_d$ be the group~ $\{{\left(
\begin{array}{cc}
    d^n \cdot \frac{s_1}{s_2} & b \\
   0 & 1 \\
    \end{array}\right)} ~|~ s_1,~s_2~ are~ nonzero~ integers~$\\
     $ coprime~ to ~d, ~b ~\in \mathbb{Q}\}$. Note $\Gamma$ is a subgroup of $G_d$ and the action of $\Gamma$ on $T_d$ can be extended to $G_d$ for  any integer prime $d >1$ (see \cite{Se}, Chaper II, Section 1.3); the case when $d$ is not a prime can be induced from the case $d$ is a prime and will be explained in the Appendix. We will assume from now on $d$ is not a prime.
On the other hand, our definition of the $\Gamma$ action on $\mathbb{R}$ can also be easily extended to $G_d$; in fact, for any $w \in \mathbb{R}$, ${\left(
\begin{array}{cc}
    d^n \cdot \frac{s_1}{s_2} & b \\
   0 & 1 \\
    \end{array}\right)}$ acts on $\mathbb{R}$ by ${\left(
\begin{array}{cc}
    d^n \cdot \frac{s_1}{s_2} & b \\
   0 & 1 \\
    \end{array}\right)} ~w = {d^n \cdot \frac{s_1}{s_2} \cdot w + b}$.  Hence we can define an action of $G_d$ on $T_d \times \mathbb{R}$ by the diagonal action for any $d > 1$.  Furthermore,  $\Gamma$ as a subgroup of $G_d$ acts on $T_d \times \mathbb{R}$ and ${(T_d \times \mathbb{R})}/\Gamma$ is homeomorphic to our previous space X. Note that this action is the ``boundary" of the natural (isometric) diagonal action of $G_d$ on $T_d \times \mathds{H}^2$ in the upper half plane model for $\mathds{H}^2$ where we identify $\mathbb{R}$ in $T_d \times \mathbb{R}$ with the $x$-axis of $\mathbb{R}^2$. Moreover, elements of the form ${\left(\begin{array}{cc}
    s & 0 \\
   0 & 1 \\
    \end{array}\right)} $ fix the line $L_0$ in $T_d$ where $s$ is an positive integer coprime to $d$, though it will not fix the whole tree $T_d$ in general; for more details about the case $d$ is a prime see \cite{Se}, Chapter II, Section 1.3, p 77, Stabilizers of straight paths: Cardan subgroups. One important observation is that the action of $g = {\left(\begin{array}{cc}
    s & 0 \\
   0 & 1 \\
    \end{array}\right)} $ on $T_d$ does not change the value of the Busemann function; i.e., $f_d(g(z)) = f_d(z)$, for every $z \in T_d$.
\end{rem}

We now define a metric on $T_d \times \mathbb{R}$ so that $\Gamma$ acts isometrically on it. Note first that both the tree $T_d$ and the real line $\mathbb{R}$ already have canonical metrics.  Now we define the metric on $T_d \times \mathbb{R}$ by the warped product of $T_d$ and $\mathbb{R}$ with respect to the warping function $d^{-f_d}$, where $f_d$ is the Busemann function we defined before. $T_d \times \mathbb{R}$ is a metric space under this metric (see for example \cite{CH} Proposition 3.1). If we restricted to the two dimensional subspace $L_0 \times \mathbb{R}$, the metric will be
$$dz^2 ~+~ (d^{-f_d(z)})^2~dw^2, ~where~z~\in ~L_0,~ w ~\in~ \mathbb{R}.$$
Simple calculations show that this is a space with constant curvature $-(\ln d)^2$. Hence $T_d \times \mathbb{R}$ is constructed by gluing together infinitely many copies of hyperbolic planes with curvature $-(\ln d)^2$ along $[z,\omega) \times \mathbb{R}$, where $z$ can be any vertex in $T_d$. In fact, each line $L$ in $T_d$ going towards $\omega$ determined a hyperbolic plane $H_L = L \times \mathbb{R}$. And if  $z \in T_d$ is the first point $L$ and $L'$ meet towards $\omega$, then $H_L$ and $H_{L'}$ are glued together along $[z,\omega)\times \mathbb{R}$. Note if we identify $L \times \mathbb{R}$ with the Poincar\'e disk model (with curvature $-(\ln d)^2$), $z \times \mathbb{R}$ will be
mapped to a horosphere. Correspondingly, $[z,\omega)$ will be mapped to the region outside the horoball (the disk
bounded by the horosphere) in the Poincar\'e disk, which is not convex. Hence, $T_d \times \mathbb{R}$ might not be a CAT(0) space, in fact,
 our group is not a CAT(0) group since it contains a non-finitely generated Abelian subgroup; cf.  Corollary 7.6 in \cite{BH}, p247.

We are going to need the following lemmas in the future.
\begin{lem}{\label{dd}}
For any two points $(z_1,w_1),(z_2,w_2) \in T_d \times \mathbb{R}$,
$$d_{T_d \times \mathbb{R}}((z_1,w_1),(z_2,w_2)) \geq d_{T_d \times \mathbb{R}}((z_1,w_2),(z_2,w_2)) = d_T(z_1,z_2) $$
\end{lem}
\Proof By Lemma 3.1 in \cite{CH}, we have $d_{T_d \times \mathbb{R}}((z_1,w_1),(z_2,w_2)) \geq d_T(z_1,z_2)$. By Lemma 3.2 in \cite{CH}, $d_{T_d \times \mathbb{R}}((z_1,w_2),(z_2,w_2)) = d_T(z_1,z_2)$. \qed

\begin{rem}
It is not always true that $d_{T_d \times \mathbb{R}}((z_1,w_1),(z_2,w_2)) \geq |w_1 - w_2|$.
\end{rem}

\begin{cor}{\label{ddc}}
Let $z_1,z_2 \in T_d$, $w_1,w_2 \in \mathbb{R}$, then the following inequality holds
$$ d_{T_d \times \mathbb{R}}((z_1,w_1),(z_2,w_2)) \geq \frac{1}{2} d_{T_d \times \mathbb{R}}((z_1,w_1),(z_1,w_2))  $$
\end{cor}
\Proof By the triangle inequality on metric spaces, we have \\
 $\hspace*{5mm} d_{T_d \times \mathbb{R}}((z_1,w_1),(z_2,w_2)) + d_{T_d \times \mathbb{R}}((z_2,w_2),(z_1,w_2)) \geq  d_{T_d \times \mathbb{R}}((z_1,w_1),(z_1,w_2))$.\\
By Lemma \ref{dd}, $d_{T_d \times \mathbb{R}}((z_1,w_1),(z_2,w_2)) \geq d_T(z_1,z_2) = d_{T_d \times \mathbb{R}}((z_2,w_2),(z_1,w_2)) $, therefore $ 2 d_{T_d \times \mathbb{R}}((z_1,w_1),(z_2,w_2)) \geq  d_{T_d \times \mathbb{R}}((z_1,w_1),(z_1,w_2))$. \qed

\begin{lem}{\label{led}}
Let $z_0$ be a fixed point in $T_d$, $w_1,w_2,w_3$ are three different points in $\mathbb{R}$ such that $|w_1 - w_2| < |w_1 - w_3|$, then
$$d_{T_d \times \mathbb{R}}((z_0,w_1),(z_0,w_2)) < d_{T_d \times \mathbb{R}}((z_0,w_1),(z_0,w_3))$$
\end{lem}

\Proof It is not hard to see that we can arrange the geodesics connecting $(z_0,w_1)$ and $(z_0,w_i)$, i =2,3, to lie within a single hyperbolic plane. Hence we only need to prove the lemma in the hyperbolic plane $L_0 \times \mathbb{R}$. We can map $L_0 \times \mathbb{R}$ to the Poincar\'e disk model with $(z_0,w_1)$ as the center; $\{z_0\} \times \mathbb{R}$ will be mapped to a horosphere. Now our lemma follows easily.
\qed

\begin{lem}\label{vl}
Let $z_0$ be a fixed point in $T_d$, $w_1,w_2$ are two fixed points in $\mathbb{R}$, denote the distance $d_{T_d \times \mathbb{R}}((z_0, \frac{w_1}{n}),(z_0, \frac{w_2}{n})) $ by $D_n$, then for any given integer $n>0$,
$$\sinh{(\frac{\ln d} {2}D_1)} = n ~\sinh{(\frac{\ln d} {2} D_n )}  $$
Therefore,
$$D_n = \frac{2}{\ln d} arcsinh(\frac{1}{n} \sinh{(\frac{\ln d} {2}D_1)})$$
In particular,
 $$\lim_{n \rightarrow \infty} D_n = 0.$$
\end{lem}
\Proof
Denote the induced inner metric on $\{z_0\} \times \mathbb{R}$ as $d$, then $d(w_1,w_2) = n ~d(\frac{w_1}{n},\frac{w_2}{n})$. Note if we use the Poincar\'e disk model, $\{z_0\} \times \mathbb{R}$ will be mapped to a horosphere, hence for any two points $w_1, w_2 \in \{z_0\} \times \mathbb{R}$,
 $$d(w_1,w_2) = \frac{2}{\ln d} \sinh{(\frac{\ln d} {2} d_{T_d \times \mathbb{R}}((z_0, w_1),(z_0, w_2)))}$$
 see for example \cite{EH}, Theorem 4.6. Now apply this to both sides of the equation $d(w_1,w_2) = n ~d_{\mathbb{R}}(\frac{w_1}{n},\frac{w_2}{n})$, our lemma follows.
\qed

\section{Flow Space for $ T_d \times {\mathbb{R}}$}
In this section we define a flow space for $E({\mathbb{Z} {[{\frac{1}{d}} ]}} \rtimes_\alpha \mathbb{Z})$,
which is a horizontal subspace of the flow space defined in Bartels and L\"uck's paper \cite{BL2}. Then we prove that it
has lots of good properties, which guarantees that it has a long thin cover.

We first introduce Bartels and L\"uck's flow space starting with the notion of generalized geodesic.
\begin{defn}
Let X be a metric space. A continuous map $c : \mathbb{R} \rightarrow X$ is called a \textbf{generalized geodesic} if
there are $c_{-}, c_{+} \in \bar{\mathbb{R}} :=  \mathbb{R} \coprod \{-\infty, \infty \}$ satisfying

$$  c_- \leq c_+, c_- \neq \infty, c_+ \neq -\infty  $$
such that $c$ is locally constant on the complement of the interval $I_c := (c_-, c_+)$ and
restricts to an isometry on $I_c$.
\end {defn}

\begin{defn}
Let $(X, d_X)$ be a metric space . Let \textbf{$FS(X)$} be the set of all generalized geodesics in
$X$. We define a metric on $FS(X)$ by

$$    d_{FS(X)}(c,d) := \int_{\mathbb{R}} \frac{d_X(c(t), d(t))}{2e^{|t|}} dt             $$

\end{defn}

The flow on $FS(X)$ is defined by
$$ \Phi: FS(X) \times \mathbb{R} \rightarrow FS(X) $$
where $ {\Phi}_{\tau} (c)(t) = c(t + \tau)$ for $ \tau \in \mathbb{R}$, $c \in FS(X)$ and $t \in \mathbb{R}$.

\begin{lem}{\label{inq}}
The map $\Phi$ is a continuous flow and if we let $c,d \in FS(X)$, $\tau \in \mathbb{R}$,
then the following inequality holds
$$ e^{-|\tau|} d_{FS(X)}(c,d) \leq d_{FS(X)}(\Phi_{\tau}(c),\Phi_{\tau}(d)) \leq e^{|\tau|} d_{FS(X)}(c,d)  $$
\end{lem}
\Proof A more general version is proved  in  \cite{BL2},  Lemma 1.3. \qed

Note that the isometry group of $(X, d_X)$ acts canonically on $FS(X)$. Recall a  map is proper if the inverse image of every compact subset is compact. Bartels and L\"uck  also proved the following for the flow space $FS(X)$ in \cite{BL2} Proposition 1.9 and 1.11.

\begin{prop} \label{ac}
If $(X,d_X)$ is a proper metric space, then $(FS(X),d_{FS(X)})$ is a proper metric space, in particular it is a complete metric space. Furthermore, if a group $\Gamma$ acts isometrically and properly on $(X, d_X)$, then $\Gamma$ also acts on $(FS(X),d_{FS(X)})$ isometrically and properly. In addition, if $\Gamma$ acts cocompactly on $X$, then $\Gamma$ acts cocompactly on $FS(X)$.

\end{prop}

Now we define our flow space by
 $$HFS(T_d \times\mathbb {R}) := FS(T_d)\times \mathbb {R}$$
  where $T_d$ has its natural metric with edge length $1$.   Since $\Gamma$ has an action on both $FS(T_d)$ and $\mathbb{R}$,
 $\Gamma$ will have a diagonal action on $FS(T_d)\times \mathbb {R}$ also. One can think of $HFS(T_d \times\mathbb {R})$ as the horizontal subspace of $FS(T_d \times\mathbb {R})$. In fact, there is a natural  embedding of  $HFS(T_d \times\mathbb {R})$ (as a topological space with product topology) into $FS(T_d \times \mathbb {R})$ defined as follows:
for a generalized geodesic $c$ on $T_d$, and $w \in \mathbb{R}$, we define a generalized geodesic on $T_d \times \mathbb{R}$, which maps $t \in \mathbb{R}$ to $(c(t),w) \in T_d \times \mathbb {R}$. $HFS(T_d \times\mathbb {R})$ will inherit a metric from this embedding.

For the rest of this section, let $X = T_d \times \mathbb {R}$.

\begin{lem}\label{prop}
The flow space $HFS(T_d \times\mathbb {R})$ is a proper metric space, in particular a complete metric space.
\end{lem}

\Proof
In order to prove $HFS(T_d \times\mathbb {R})$ is a proper metric space, we need to show every closed ball $B_r(c) ~= ~\{c' ~|~ d_{HFS(T_d \times\mathbb {R})}(c,c') \leq r\}$ in $HFS(T_d \times\mathbb {R})$ is compact. Let $\{c_i\}$ be a Cauchy sequence in the closed ball $B_r(c)$, we need to show it converges to a point in $B_r(c)$. Since the space $FS(T_d \times \mathbb {R})$ is proper, we can now assume $\{c_i\}$ converges to a point $c_0$ in $FS(T_d \times \mathbb {R})$. We only need to show $c_0 \in HFS(T_d \times\mathbb {R})$. Denote the projection map from $T_d \times \mathbb{R}$ to $T_d$ as $q_1$, from $T_d \times \mathbb{R}$ to $\mathbb{R}$ as $q_2$, then $c_i(t) = (q_1(c_i(t),q_2(c_i(t)))$. Suppose $c_0 \notin HFS(T_d \times\mathbb {R})$, then $q_2(c_0(t))$ is not a constant map. Choose a big enough close interval $I$ in $\mathbb{R}$ such that $q_2(c_0(t))$ restricted to $I$ is not a constant, we can assume the maximum value is $A_1$, while the minimum is $A_2$, where $A_1 > A_2$. Let $\delta = A_1 - A_2 >0$ and $I_1 = \{t \in I~|~ q_2(c_0(t)) > A_1 - \frac{\delta}{4}\}$, correspondingly $I_2 = \{t \in I~|~ q_2(c_0(t)) < A_2 + \frac{\delta}{4}\}$. Note $I_1$ and $I_2$ are nonempty sets with measure bigger than $0$. Now for any given $c_i$, if $q_2(c_i) \geq \frac{A_1 + A_2}{2}$, then

$$d_{FS(X)}(c_0, c_i) = \int_{\mathbb{R}} \frac{d_{T_d \times \mathbb{R}}(c_0(t), c_i(t))}{2e^{|t|}} dt \hspace*{69mm}$$

$$\geq \int_{I_2} \frac{d_{T_d \times \mathbb{R}}(c_0(t), c_i(t))}{2e^{|t|}} dt \hspace*{46mm}$$
$$\geq \int_{I_2} \frac{d_{T_d \times \mathbb{R}}((q_1(c_0(t)),q_2(c_0(t))), (q_1(c_i(t)),q_2(c_i(t))))}{2e^{|t|}} dt $$
$$\hspace*{13mm} = \int_{I_2} \frac{\frac{1}{2}d_{T_d \times \mathbb{R}}(c_0(t), (q_1(c_0(t)), q_2(c_i(t))))}{2e^{|t|}} dt\hspace*{5mm}( by~ Corollary~ \ref{ddc})$$
$$\hspace*{13mm} \geq \int_{I_2} \frac{\frac{1}{2}d_{T_d \times \mathbb{R}}(c_0(t), (q_1(c_0(t)), \frac{A_1 +A_2}{2})))}{2e^{|t|}} dt > 0\hspace*{5mm}( by~ Lemma~ \ref{led})$$
The last integral is independent of $c_i$; denote its value as $\epsilon_1$. For the same reason, if $q_2(c_i) < \frac{A_1 + A_2}{2}$, there exists  $\epsilon_2 > 0$ such that $d_{FS(X)}(c_0, c_i) \geq \epsilon_2$. Let $\epsilon = \min(\epsilon_1,\epsilon_2) > 0$, then $d_{FS(X)}(c_0, c_i) \geq \epsilon >0 $.
Hence the sequence $\{c_i\}$ can never converge to $c_0$, contradiction. \qed

\begin{rem}
The proof in fact shows that the embedding $HFS(X) \subset FS(X)$ is a closed $\Gamma$-equivariant embedding.
\end{rem}

We define now the flow
$$ \Phi: HFS(T_d \times\mathbb{R}) \times \mathbb{R} \rightarrow HFS(T_d \times\mathbb{R}) $$
by $ {\Phi}_{\tau} ((c,w))(t) = (c(t + \tau),w)$ for $c \in FS(T_d)$ and $ \tau,w,t \in \mathbb{R}$. Note $\Phi$ is a $\Gamma$-equivariant flow.
~\\

\begin{lem} {\label{l3}}
The flow space $HFS(X)$ has the following properties:
\begin{itemize}
    \item[(i)] $\Gamma$ acts properly and cocompactly on $HFS(X)$.
    \item [(ii)]  Given $C > 0$, there are only finitely many $\Gamma$ orbits of periodic flow curves with period less than $C$ (but bigger than $0$).
    \item [(iii)] Let $HFS(X)^{\mathbb{R}}$ denote the $\mathbb{R}$-fixed point set, i.e., the set of points $c \in HFS(X)$ for which ${\Phi}_{\tau} (c) = c$ for all $\tau \in \mathbb{R}$, then $HFS(X) - HFS(X)^{\mathbb{R}}$ is locally connected.
    \item[(iv)] If we put
    $$ k_\Gamma := sup\{ |H| | H \subseteq \Gamma subgroup~with ~finite~ order~ |H| \} ;$$
    $$d_{HFS(X)} := dim(HFS(X) - {HFS(X)}^{\mathbb{R}}); \hspace*{37mm}$$
    then $k_\Gamma$ and $d_{HFS(X)}$ are finite.

\end{itemize}

\end{lem}

\Proof Those properties are essentially implied by results in \cite{BL2}, section 1.

\begin{itemize}
 \item[(i)]  The $\Gamma$ action on $HFS(X)$ is proper since $HFS(X)$ $\Gamma$-equivariantly embeds into $FS(X)$ and the $\Gamma$-action on $FS(X)$ is proper by Proposition \ref{ac}. By Lemma 1.10 in \cite{BL2}, for any $w_0 \in \mathbb{R}$, the evaluation map $FS(T_d) \rightarrow T_d$ defined by $c \rightarrow c(w_0)$ is proper. Hence $FS(T_d) \times \mathbb{R} \rightarrow T_d \times \mathbb{R}$ is proper. This will induce a map $(FS(T_d) \times \mathbb{R})/{\Gamma} \rightarrow (T_d \times \mathbb{R})/{\Gamma}$. Since $(T_d \times \mathbb{R})/{\Gamma}$ is compact, $(FS(T_d) \times \mathbb{R})/{\Gamma}$ is compact also. One can prove this by choosing a compact fundamental domain in $T_d \times \mathbb{R}$, and using the fact that the map $FS(T_d) \times \mathbb{R} \rightarrow T_d \times \mathbb{R}$ is proper. Hence $\Gamma$ acts cocompactly on $HFS(X)$.

\item [(ii)] Note that periodic orbits in $HFS(T_d\times \mathbb{R})$ are periodic orbits in $FS(T_d\times \mathbb{R})$, which move horizontally (i.e., move along the tree direction, with $\mathbb{R}$ coordinate fixed). Note also that the embedding of $HFS(T_d\times \mathbb{R})$ into $FS(T_d\times \mathbb{R})$ is a $\Gamma$-equivariant map, and there are only finitely many nonzero horizontal periodical geodesics on $(T_d\times \mathbb{R})/{\Gamma}$ of period less than $C$. In fact $(T_d\times\mathbb{R})/{\Gamma} = S^1 \times [0,1] / (z,0) \thicksim (z^d, 1)$ (see Section \ref{model}) and horizontal periodical geodesics with period $m$ on it corresponding to solutions of the equation $d^m x ~\equiv~ x~ (~mod~ 1~)$, where $x \in \mathbb{R}/\mathbb{Z}$. Since for any positive integer $m$, the equation has finitely many solutions. In fact the solutions are $x = \frac{k}{d^m- 1}, ~k \in \mathbb{Z}$, $x\in\mathbb{R}/\mathbb{Z}$. Hence for fixed $m$, the number of solutions is $d^m- 1$. So there will be only finitely many nonzero horizontal periodical orbits on $HFS((T_d\times\mathbb{R})/{\Gamma} )$ with period less than $C$. Hence the claim in (ii) now follows.

 \item [(iii)] This is because $T_d$ is a tree, hence a CAT(0) space. In Bartels and L\"uck's paper \cite{BL2}, Proposition 2.10, they proved that for any CAT(0) space A, $FS(A) - FS(A)^{\mathbb{R}}$ is locally connected. Our flow space $HFS(X) = FS(T_d) \times \mathbb{R}$. Since the flow on our flow space only flow on the first factor, $HFS(X)^{\mathbb{R}} = FS(T_d)^{\mathbb{R}} \times \mathbb{R}$. So $HFS(X) - HFS(X)^{\mathbb{R}}$ will be locally connected as well.

 \item[(iv)] Since $\Gamma$ is a torsion free group, $k_\Gamma = 1$. Note that any two point can be connected by a unique geodesic in the tree $T_d$, it is not hard to see that the flow space $FS(T_d)$ will have dimension  less than $5$. Therefore our flow space $HFS(X) = FS(T_d) \times \mathbb{R}$ has finite dimension, and     hence  $d_{HFS(X)}$ is finite. One can also consult Bartels and L\"uck's paper \cite{BL2}, they have more general results for CAT(0) spaces.
 \end{itemize}
\qed

\begin{rem} \label{remfl}
We define an embedding $\Psi: T_d \times \mathbb{R} \rightarrow FS(T_d) \times \mathbb{R}$, by $(z,w) \rightarrow (c_z,w)$, where $c_z$
is the unique generalized geodesic which sends $(-\infty, 0)$ to $z$, and $[0,\infty)$ isometrically to the geodesic $[z,\omega)$. Recall  $[z,\omega)$ is the unique geodesic connecting $z$ and the specified infinity point $\omega$ of $T_d$ defined in Section \ref{model}. Also, we can flow this embedding by flowing its image in $HFS(X)$; define $\Psi_{\tau}(z,w) = \Phi_{\tau}(\Psi(z,w))$. It is easy to see that $\Psi_{\tau}$ is a $\Gamma$-equivariant map since $\omega$ is fixed under the group action.
\end{rem}

We need the following lemma in the future.
\begin{lem} \label{lpar}
Let $z_0$ be a fixed point in $T_d$, $w_1,w_2$ are two fixed points in $\mathbb{R}$, and $P_n = (z_0,\frac{w_1}{n})$, $Q_n = (z_0,\frac{w_2}{n})$, $d_X(P_1,Q_1) < D$. Then for any $\epsilon > 0$, there exists a number $\bar{N}$, which depends only on $\epsilon$, $D$ and $d$, such that for any $n > \bar{N}$
$$d_X(P_n,Q_n) < \frac{\epsilon}{4}$$
and
$$ d_{HFS(X)}(\Psi(P_n),\Psi(Q_n)) \leq \epsilon$$
\end{lem}
\Proof Choose $T= \max\{1,\ln {\frac{4}{\epsilon}}\}$ . Since $\Psi(P_n)(T)$ and $\Psi(Q_n)(T)$ have the same $T_d$ coordinate, by Lemma \ref{vl}, we can choose a big enough integer $N$ such that for any $n > N$, $d_X(\Psi(P_n)(T),\Psi(Q_n)(T)) < \frac{\epsilon}{4}$. Note $\bar{N}$ depends only on $\epsilon$, $D$ and $d$.  Using the definition of generalized geodesic, we have $d_{X}( \Psi(P_n)(t), \Psi(P_n)(T)) = t-T$ and $d_{X}(\Psi(Q_n)(t), \Psi(Q_n)(T)) = t-T$, where $t \geq T$. Hence for any $t \geq T$, by triangle inequality,  we have the following
$$d_{X}(\Psi(P_n)(t)),\Psi(Q_n)(t))\hspace*{88mm} $$
$$\leq d_{X}( \Psi(P_n)(t),  \Psi(P_n)(T)) + d_{X}(\Psi(P_n)(T),\Psi(Q_n)(T))  + d_{X}(\Psi(Q_n)(T), \Psi(Q_n)(t)) $$
$$\leq \frac{\epsilon}{4} + 2(t-T) \hspace*{103mm}$$
On the other hand, the metric defined on $T_d \times \mathbb{R}$ is expanding in the $\mathbb{R}$ direction when moving towards $\omega$. Hence for any $0\leq t \leq T$, we have
$$d_{X}(\Psi(P_n)(t)),\Psi(Q_n)(t)) \leq d_{X}(\Psi(P_n)(T)),\Psi(Q_n)(T)) < \frac{\epsilon}{4}$$
And for $t \leq 0$, $\Psi(P_n)(t) = \Psi(P_n)(0) = P_n$, $\Psi(Q_n)(t) = \Psi(Q_n)(0) = Q_n$, hence
$$d_{X}(\Psi(P_n)(t)),\Psi(Q_n)(t)) =  d_{X}(P_n,Q_n) < \frac{\epsilon}{4}.$$

 Therefore, for any $n>\bar{N}$
$$d_{HFS(X)}(\Psi(P_n),\Psi(Q_n)) = \int_{\mathbb{R}} \frac{d_X(\Psi(P_n)(t)),\Psi(Q_n)(t)))}{2e^{|t|}} dt \hspace*{30mm} $$
$$ = \int_{(-\infty,0]} \frac{d_X(\Psi(P_n)(t)),\Psi(Q_n)(t)))}{2e^{|t|}} dt ~~+~~ \int_{[0,T]} \frac{d_X(\Psi(P_n)(t)),\Psi(Q_n)(t)))}{2e^{|t|}} dt$$

$$\hspace*{20mm}+ \int_{[T, \infty)} \frac{d_X(\Psi(P_n)(t)),\Psi(Q_n)(t)))}{2e^{|t|}} dt$$

$$ \leq \int_{(-\infty,0]} \frac{\frac{\epsilon}{4}}{2e^{|t|}} dt \hspace*{2mm}+ \hspace*{2mm}\int_{[0,T]} \frac{\frac{\epsilon}{4}}{2e^{|t|}} dt  \hspace*{2mm}+\hspace*{2mm} \int_{[T,\infty)} \frac{\frac{\epsilon}{4} + 2(t-T)}{2e^{|t|}} dt \hspace*{20mm}$$
$$= \frac{\epsilon}{4} +  e^{-T} \leq \frac{\epsilon}{4} + \frac{\epsilon}{4} = \frac{\epsilon}{2}\hspace*{75mm}$$
Hence we proved the Lemma.
\qed

Because of the properties proved in Lemma \ref{l3}, Theorem 1.4 in \cite{BLR} yields  a long thin cover for $HFS(X)$; i.e. the following result holds

\begin{prop}\label{ltc}
 There exists a natural number $N$, depending only on $k_\Gamma,~d_{HFS(X)}$ and the action of $\Gamma$ on an arbitrary neighborhood of ${HFS(X)}^{\mathbb{R}}$ such that for every $\lambda > 0$ there is an $\mathcal{VC}yc$-cover $\mathcal{U}$ of $HFS(X)$ with the following properties:
 \begin{itemize}
  \item[(i)] dim $\mathcal{U} \leq N$;
 \item [(ii)] For every $x \in HFS(X)$ there exists $U_x \in \mathcal{U}$ such that
 $$ \Phi_{[-\lambda, \lambda]}(x) := \{ \Phi_{\tau}(x) ~|~ \tau \in [-\lambda, \lambda]  \}   \subseteq U_x; $$

  \item [(iii)] ${\Gamma} \setminus {\mathcal{U}}$ is finite.
 \end{itemize}
where $\mathcal{VC}yc$ denote the collections of virtually cyclic subgroups of a group.
 \end{prop}



Recall that the dimension of a cover $\mathcal{U}$ is defined to be the greatest $N$ such that there exists $N+1$ elements in $\mathcal{U}$ with nonempty intersection.  In general, for a collection of subgroups $\mathcal{F}$, we define a $\mathcal{F}$-cover as following.

\begin{defn}
Let G be a group and Z be a G-space. Let $\mathcal{F}$ be a collection of subgroups of G. An open cover $\mathcal{U}$ of Z is called an \textbf{$\mathcal{F}$-cover} if the following three conditions are satisfied.\\
(i) For $g\in G$ and $ U \in \mathcal{U}$ we have either $g(U) = U$ or $g(U) \bigcap U = \emptyset$;\\
(ii) For $g\in G$ and $ U \in \mathcal{U}$, we have $g(U) \in \mathcal{U}$;\\
(iii) For $U \in \mathcal{U}$ the subgroup $G_U := \{ g \in G~|~ g(U) = U \}$ is a member of $\mathcal{F}$.

\end{defn}

For a subset A of a metric space Z and $\delta >0$, $A^{\delta}$ denotes the set of all points $z \in Z$ for which $d(z,A) < \delta$. Combining Lemma \ref{inq} and the fact that $\Gamma$ acts cocompactly on ${FS(T_d) \times \mathbb{R}}$ (Lemma \ref{l3}, (i)), Proposition \ref{ltc} can be improved to the following.

\begin{prop}\label{lte}
 There exists a natural number $N$, depending only on $k_\Gamma,d_{HFS(X)}$ and the action of $\Gamma$ on an arbitrary neighborhood of ${HFS(X)}^{\mathbb{R}}$ such that for every $\lambda > 0$ there is a $\mathcal{VC}yc$-cover $\mathcal{U}$ of $HFS(X)$ with the following properties:
 \begin{itemize}
  \item[(i)] dim $\mathcal{U} \leq N$;
 \item [(ii)] There exists a $\delta >0$ depends on $\lambda$ such that for every $x \in HFS(X)$ there exists $U_x \in \mathcal{U}$ such that
 $$ (\Phi_{[-\lambda, \lambda]}(x))^{\delta}   ~\subseteq~  U_x; ~~~~~~~$$

  \item [(iii)] ${\Gamma} \setminus {\mathcal{U}}$ is finite.
 \end{itemize}
 \end{prop}
\Proof A simple modification of the argument in \cite{BLR}, section 1.3, page 1804-1805, yields the result. In their proof, they used a lemma (Lemma 7.2) which will be replaced by Lemma \ref{inq} in our case. \qed


\section {Hyper-elementary subgroups of ${\mathbb{Z}}_{q^s} \rtimes_\alpha {\mathbb{Z}}_{t_s}$} {\label{gts}}
In this section, we study the hyper-elementary subgroups of ${\mathbb{Z}}_{q^s} \rtimes_\alpha {\mathbb{Z}}_{t_s}$, where $d$ is an integer such that $|d| >1$,  $q$ is a prime number greater than $|d| +1$, $s$ is a positive integer, $\alpha$ is multiplication by $d$, $t_s$ is the order of $d$ in the group of units of the ring ${\mathbb{Z}}_{q^s}$. The reason we study this is ${\mathbb{Z}}_{q^s} \rtimes_\alpha {\mathbb{Z}}_{t_s}$ can be realized as a quotient of  ${\mathbb{Z} {[{\frac{1}{d}} ]}} \rtimes_d \mathbb{Z}$. In fact, note first $(q^s {\mathbb{Z} {[{\frac{1}{d}} ]}}) \rtimes_d \{0\}$ is a normal subgroup of ${\mathbb{Z} {[{\frac{1}{d}} ]}} \rtimes_d \mathbb{Z}$, therefore after we modulo it out, the quotient group is ${\mathbb{Z}}_{q^s} \rtimes_d {\mathbb{Z}}$. $1 \in \mathbb{Z}$ acts on ${\mathbb{Z}}_{q^s}$ as multiplication by $d$, which has order $t_s$ when considered as an element of $Aut({\mathbb{Z}}_{q^s})$. Hence we can further map ${\mathbb{Z}}_{q^s} \rtimes_d {\mathbb{Z}}$ to  ${\mathbb{Z}}_{q^s} \rtimes_d {\mathbb{Z}}_{t_s}$.

We denote the group of units of a ring $R$ by $U(R)$, hence $t_s$ is the order of $d$ in $U({\mathbb{Z}}_{q^s})$.

\begin{defn}
A \textbf{hyper-elementary group} H is an extension of a $p$-group by a cyclic group of order n, where $p$ is a prime number, $(n, p) = 1$, in other words,
there exists a short exact sequence
$$ 1 \rightarrow C_n \rightarrow H \rightarrow G_p \rightarrow 1 $$
where $C_n$ is a cyclic group of order n, $G_p$ is a $p$-group such that $(n, p) = 1$.
\end{defn}

Note first that the order of the units group $U({\mathbb{Z}}_{q})$
 is $q - 1$,  while the order of $U({\mathbb{Z}}_{q^s})$ is $q^s - q^{s - 1} = q^{s - 1}(q - 1)$.

\begin{lem} \label{l1}

If $d^t \equiv 1 ~(~mod~ q~)$, then $d^{tq^{s-1}} \equiv 1 (~mod~ q^s~)$.

\end{lem}

\Proof We prove this by induction. For $s = 1$, this is automatically true. Now assume it is true for $k$, we prove it for $k + 1$.
By hypothesis, $d^{tq^{k-1}} = mq^k + 1$,where $m$ is an integer. So $d^{tq^k} = (mq^k + 1)^q$. Expanding the right side, it is easy to see that
$d^{tq^{k}} \equiv 1 (~mod~ q^{k +1}~) $.
\qed
\\

If we assume the order of $d$ in $U({\mathbb{Z}}_{q})$ is $t_1$, and the order of $d$ in $U({\mathbb{Z}}_{q^s})$ is $t_s$, then $t_1 ~|~ q - 1$, and $t_s ~|~ q^{s - 1}(q - 1)$. Let $t_s = {m_s}q^{k_s}$, where $(m_s, q) = 1$. Then what the previous lemma says is that $m_s ~|~ t_1$ and $k_s \leq s - 1$. In the following lemma, we prove $m_s = t_1$.

\begin{lem}\label{l2}
Assume $t_1$ is the order of $d$ in $U({\mathbb{Z}}_{q})$, $t_s = {m_s}q^{k_s}$ is the order of $d$ in $U({\mathbb{Z}}_{q^s})$, where $(m_s, q) = 1$, then $m_s = t_1$.
\end{lem}

\Proof We prove it by contradiction. Assume $m_s \neq t_1$, denote $d^{m_s} = a$, then $a$ is not equal to $1$ modulo $q$, since $m_s ~|~ t_1$. But since the order of $d$ in
$U({\mathbb{Z}}_{q^s})$ is ${m_s}q^{k_s}$, we have ~$d^{{m_s}q^{k_s}} \equiv~  a^{q^{k_s}} \equiv 1 ~(~mod ~q^s~)$. This means $a^{q^{k_s}} \equiv 1 ~(~mod ~q~)$, which is not true. In fact, by Fermat's theorem, $a^{q^{ k_s}} \equiv a^{q^{k_s - 1}} \equiv \ldots a ~(~mod ~q~)$.
\qed

We list some basic formulas we are going to use frequently:\\
(1) $(a,b)(a',b') = (a+{d^{-b}}a', b+b')$;~\\
(2) ${(a,b)}^n = ((1 + d^{-b} + \ldots + d^{-(n - 1)b})a,nb)$;
~\\
(3) ${(a,b)}^{-1} = (-{d^b}a, -b)$;
~\\
(4) $(x,y)(a,b){(x,y)}^{-1} = ((1 - d^{-b})x + {d^{-y}}a,b)$;\\
where $(a,b),(a',b'),(x,y) \in {\mathbb{Z}}_{q^s} \rtimes_\alpha {\mathbb{Z}}_{t_s}$,  $n$ is a positive integer and $d^{-b}$ is the inverse element of $d^b$ in the unit group $U({\mathbb{Z}}_{q^s})$, etc.

\begin{lem} \label{conj}
For any $(a,b) \in {\mathbb{Z}}_{q^s} \rtimes_\alpha {\mathbb{Z}}_{t_s}$, if $t_1 \nmid b$, then $(a,b)$ can be conjugate to $(0,b)$.
\end{lem}
\Proof
 Note first if $t_1 \nmid b$, then $d^{-b} \nequiv 1 ~(~mod~ q~) $, hence $1 - d^{-b}$ is a unit in $Z_{q^s}$. Since $(x,y)(a,b){(x,y)}^{-1} = ((1 - d^{-b})x + {d^{-y}}a,b)$, if $1 - d^{-b}$ is a unit in $Z_{q^s}$, then we can find $x,y$, such that $(1 - d^{-b})x + {d^{-y}}a ~=~0$, which means we can conjugate  $(a,b)$ to $(0,b)$. In fact, we can always take $y = 0$.
\qed

\begin{thm}\label{thm3}
For $s>q$ and $q$ a prime number bigger than $|d|+1$, every hyper-elementary subgroup H of ${\mathbb{Z}}_{q^s} \rtimes_\alpha {\mathbb{Z}}_{t_s}$ can be conjugated to a subgroup of one of the following three types of subgroups: \\
\hspace*{5mm} type 1:\hspace*{2mm} ${\mathbb{Z}}_{q^s} \rtimes_\alpha {\mathbb{Z}}_{q^{k_s}}$;\\
\hspace*{5mm} type 2:\hspace*{2mm} ${\mathbb{Z}}_{q^s} \rtimes_\alpha {\mathbb{Z}}_{t_1}$;\\
\hspace*{5mm} type 3:\hspace*{2mm} $\{0\}\rtimes_\alpha{\mathbb{Z}}_{t_s}$.\\
where  $\alpha$ is multiplication by $d$, $t_1$ is the order of $d$ in the group of units $U({\mathbb{Z}}_q)$ and $t_s$ is the order of $d$ in the group of units $U({\mathbb{Z}}_{q^s})$, $t_s = t_1 q^{k_s}$.
\end{thm}

\begin{rem}
The homomorphism $\alpha: {\mathbb{Z}}_{q^{k_s}} \rightarrow Aut({\mathbb{Z}}_{q^s})$ in  ${\mathbb{Z}}_{q^s} \rtimes_\alpha {\mathbb{Z}}_{q^{k_s}}$ is the restriction of the homomorphism $\alpha: {\mathbb{Z}}_{t_s} \rightarrow Aut({\mathbb{Z}}_{q^s})$, to the subgroup
${\mathbb{Z}}_{q^{k_s}} \subseteq {\mathbb{Z}}_{t_s}$. There is a similar remark for ${\mathbb{Z}}_{q^s} \rtimes_\alpha {\mathbb{Z}}_{t_1}$.
\end{rem}

\Proof Denote the cyclic part of the hyper-elementary subgroup H by C. If $s > q$, then $q^s > |d|^{q-1}$, which implies $k_s \geq 1$.

First, if C is a trivial group, then H is just a $p$-group. By Sylow's theorem, any $p$-group can be conjugate to a subgroup of a maximal $p$-group. Note that the order of ${\mathbb{Z}}_{q^s} \rtimes_\alpha {\mathbb{Z}}_{t_s}$ is $q^s \cdot t_s= q^s\cdot {t_1}q^{k_s} = t_1q^{s+k_s}$, where $t_1$ and $q$ are coprime, hence any $p$-group in ${\mathbb{Z}}_{q^s} \rtimes_\alpha {\mathbb{Z}}_{t_s}$ can be conjugated to a subgroup of ${\mathbb{Z}}_{q^s} \rtimes_\alpha {\mathbb{Z}}_{q^{k_s}}$ (type 1), or $\{0\}\rtimes_\alpha{\mathbb{Z}}_{t_1}$ which is a subgroup of $\{0\}\rtimes_\alpha{\mathbb{Z}}_{t_s}$ (type 3).

 Now assume C is non-trivial, say C is generated by $(a,b) \in {\mathbb{Z}}_{q^s} \rtimes_\alpha {\mathbb{Z}}_{t_s}$.

 If $b = 0$, C will have order a power of $q$, then since the $p$-group belonging to $H$ has to be coprime to q, and its order has to divide $t_1$, which means H is a subgroup of ${\mathbb{Z}}_{q^s} \rtimes_\alpha {\mathbb{Z}}_{t_1}$ (type 2). So for the rest of the proof, we will assume that $b \neq 0$.

 If $a = 0$,  C lies in $\{0\} \rtimes_\alpha {\mathbb{Z}}_{t_s}$. If the $p$-group part of H is trivial, then we are in type 3. Otherwise in order for C to be a normal subgroup, for any $(x,y) \in H$, $(x,y)(0,b){(x,y)}^{-1} = ((1 - d^{-b})x, b)$ has to lie in C as well, which means $(1 - d^{-b})x$ has to be equal to zero in ${\mathbb{Z}}_{q^s}$. Note first that this means C lies in the center of H. If $x = 0$ for all $(x,y)$, then H is a subgroup of $\{0\} \rtimes_\alpha {\mathbb{Z}}_{t^s} $ (type 3). If $x \neq 0$ for some $(x,y)$ which lies in the $p$-group of H, then at least, $(1 - d^{-b}) \equiv 0 ~(~mod~ q~)$ and $x  \equiv 0 ~(~mod~ q~)$. Hence $d^{-b} \equiv 1 ~(~mod~ q~)$, therefore $t_1~|~b$. So the order of the cyclic group is a power of $q$, while the order of H's $p$-group has to be coprime to $q$. By Sylow's theorem, H's $p$-group can be conjugate to a subgroup of $\{0\}\rtimes_\alpha{\mathbb{Z}}_{t_1}$, hence it is actually a cyclic subgroup. This implies the hyper-elementary subgroup H is an abelian group since C is in the center of H. Combined with the fact that the order of the cyclic group and $p$-group are coprime, we have H is again a cyclic group, say generated by $(a',b')$, where $b'$ is not zero. The case $a' = 0$ is contained in type 3. For the case $a' \neq 0$, if $t_1~\nmid~b'$, by Lemma \ref{conj}, $(a',b')$ can be conjugate to $(0,b')$ (type 3). When $t_1 ~|~ b'$, H will lie in ${\mathbb{Z}}_{q^s} \rtimes_\alpha {\mathbb{Z}}_{q^{k_s}}$ (type 1).

Now suppose both $a$ and $b$ are not zero. If $t_1$ does not divide $b$, then by Lemma \ref{conj}, we can conjugate  $(a,b)$ to $(0,b)$, which can be included in the $a = 0$ case. When $t_1~|~b$, then $d^{-b} \equiv 1 ~(~mod~ q~) $. And $b$ lies in ${\mathbb{Z}}_{q^{k_s}}$($\subseteq  {\mathbb{Z}}_{t_s}$), hence the order of C is a power of $q$ since $(a,b)$ generates C in ${\mathbb{Z}}_{q^s} \rtimes_\alpha {\mathbb{Z}}_{t_s}$.  Let $b$'s order in ${\mathbb{Z}}_{t^s}$ be $q^r$.  Since H's $p$-group has order coprime to $q$, again by Sylow's theorem, it can be conjugate to a subgroup of $\{0\}\rtimes_\alpha{\mathbb{Z}}_{t_1}$. Therefore we can assume the generator of H's $p$-group to be $(0,y)$. If $y \neq 0$, then $t_1$ does not divide $y$, so $d^{y}$ is not equal to $1$ modulo $q$. On the other hand, since $(0,y)(a,b){(0,y)}^{-1} = ({d^{-y}}a,b)$, $({d^{-y}}a,b)$ has to lie in C, which means  $({d^{-y}}a,b) = {(a,b)}^n$ for some $n$. Recall from formula (2) that  ${(a,b)}^n = ((1 + d^{-b} + \ldots + d^{-(n - 1)b})a,nb)$. In order for $nb = b$, which is the same as $(n-1)b = 0$, $b$'s order $q^r$ has to divide $n - 1$. For such $n$, $\Sigma = 1 + d^{-b} + \ldots + d^{-(n - 1)b}$  will equal to $1$ modulo $q$; note $d^{-b} \equiv 1 ~(~mod~ q~) $.  Consequently, $d^{-y} - \Sigma$ is a  unit in ${\mathbb{Z}}_{q^{k_s}}$ and $(d^{-y} - \Sigma)a = 0 \in {\mathbb{Z}}_{q^{k_s}}$, which is a contradiction since $a \neq 0$. Therefore, $y$ has to be $0$ and H's $p$-group has to be trivial, hence H will now be a subgroup of ${\mathbb{Z}}_{q^s} \rtimes_\alpha {\mathbb{Z}}_{q^{k_s}}$ (type 1). \qed

\begin{cor}
For any $n > 1$, $q$ a prime number greater than $|d|^n$ and $s > q$,  each hyper-elementary subgroup
of  ${\mathbb{Z}}_{q^s} \rtimes_\alpha {\mathbb{Z}}_{t_s}$ has index greater than $n$. Recall  $\alpha$ is multiplication by $d$ and $t_s$ is the order of $d$ in the group of units $U({\mathbb{Z}}_{q^s})$.
\end{cor}

\begin{cor} \label{c2}

For any $n > 1$, $q$ a prime number greater than $|d|^n$ and $s > q$, each hyper-elementary subgroup
of  ${\mathbb{Z}}_{q^s} \rtimes_\alpha {\mathbb{Z}}_{t_s}$ is conjugate to a subgroup H, such that one of the following is true:
\\
(1) the index $[{\mathbb{Z}}_{t_s}, \pi(H)] \geq n$, where $\pi : {\mathbb{Z}}_{q^s} \rtimes_\alpha {\mathbb{Z}}_{t_s} \mapsto {\mathbb{Z}}_{t_s}$ is the natural epimorphism.\\
(2) H is a subgroup of $\{0\}\rtimes_\alpha{\mathbb{Z}}_{t_s}$, and $q^s \geq n$.
\end{cor}

\Proof The index is $t_1$, $q^{k_s}$ for subgroups of type $1$, type $2$ respectively, and $q^s$ for subgroups of type $3$ in Theorem \ref{thm3}. If $q > |d|^n$, then the order of $d$ in ${\mathbb{Z}}_{q}$
is greater than $n$. Also $s > q$, so $k_s \geq 1$. Hence the order of $d$ in ${\mathbb{Z}}_{q^s}$ is $t_s = {t_1}q^{k_s}$, where $t_1 > n$, $q >  n$, and $k_s \geq 1$. Now by Theorem \ref{thm3} the two corollaries follow easily.
\qed

\section{Proof of the main Theorem}
In this section we prove our main theorem. Our strategy is to prove that ${\mathbb{Z} {[{\frac{1}{d}} ]}} \rtimes_\alpha \mathbb{Z}$ is in fact a Farrell-Hsiang group, as defined by  Bartels and L\"uck in \cite{BL1}. Recall that $d$ is a positive integer greater than $1$ and $\alpha$ is multiplication by $d$.

\begin{defn}
Let $\mathcal{F}$ be a family of subgroups of the finitely generated group G. We call G a \textbf{Farrell-Hsiang Group} with respect to the family $\mathcal{F}$ if the following holds for a fixed word metric $d_G$:\\
  ~~ There exists a  fixed natural number $N$ such that for every natural number $n$ there is a surjective homomorphism $\Delta_n :G \rightarrow F_n$ with $F_n$ a finite group such that the following condition is satisfied. For every hyper-elementary subgroup H of $F_n$ we set $\bar{H} := \Delta^{-1}_n (H)$ and require that there exists a simplicial complex $E_H$ of dimension at most N with a cell preserving simplicial $\bar{H}$-action whose stabilizers belong to $\mathcal{F}$, and an $\bar{H}$-equivariant map $f_H: G \rightarrow E_H$ such that $d_G(g_0,g_1) < n$ implies $d^1_{E_H}(f_H(g_0),f_H(g_1)) < \frac{1}{n}$ for all $g_0,g_1 \in G$, where $d^1_{E_H}$ is the $l^1$-metric on $E_H$.

\end{defn}
\begin{rem}
As pointed out in \cite{BFL}, Remark 1.15, in order to check a group G is a Farrell-Hsiang group, it suffices to check these conditions for one hyper-elementary subgroup in every conjugacy class of such subgroups of $F_n$.
\end{rem}

With this definition, they proved the following theorem:
\begin{thm}
Let G be a Farrell-Hsiang group with respect to the family $\mathcal{F}$. Then G satisfies the K-theoretic and L-theoretic Farrell-Jones Conjecture
with respect to the family $\mathcal{F}$.
\end{thm}

Since both the K-theoretic and L-theoretic  Farrell-Jones conjecture have been verified for abelian groups with respect to the family of virtually cyclic subgroups, by the transitivity principle (see for example \cite {BFL}, Theorem 1.11), in order to prove our main theorem, we only need to prove it with respect to the family of abelian subgroups. \\

\textbf{Claim}: The group $\Gamma = {\mathbb{Z} {[{\frac{1}{d}} ]}} \rtimes_\alpha \mathbb{Z}$ is a Farrell-Hsiang group with respect to the family of abelian subgroups.\\

 Choose $N > 0$ to be the number which appears in Proposition \ref{lte}, note $N$ is independent of $\lambda$. First, there is a quotient map ${\mathbb{Z} {[{\frac{1}{d}} ]}} \rtimes_\alpha \mathbb{Z} \rightarrow  {\mathbb{Z}}_{q^s} \rtimes_\alpha {\mathbb{Z}}_{t_s}$,  where $t_s$ is the order of $d$ in the unit group of ${\mathbb{Z}}_{q^s}$. Let $F_n$ be ${\mathbb{Z}}_{q^s} \rtimes_\alpha {\mathbb{Z}}_{t_s}$, where we choose $q$ to be a prime number bigger than $d^n$ and $s >q$ (Hence Corollary \ref{c2} holds; we will manipulate $s$ more in the future to get more control.) Let $\Delta_n$ be the quotient map. Also let H be a hyper-elementary subgroup of $F_n$ and $\bar{H} = \Delta^{-1}_n (H)$. For convenience, the metric we put on $\Gamma$ is the one inherited from the orbit embedding $\eta: \Gamma \rightarrow T_d \times \mathbb{R}$ where $\eta(g) = gx_0$, $g \in \Gamma$ and $x_0 = (P_0,0) \in T_d \times \mathbb{R} $ is the base point, denote this metric on $\Gamma$ by $d_{\Gamma}$. This metric is quasi-isometric to any word metric by Milnor-\u{S}varcz lemma (see \cite{BH},Proposition 8.19, pp140), hence it is good enough for our purpose.

By Corollary \ref{c2}, For any $n > 0$, $q$ a prime number greater than $d^{4n^2}$ and $s > q$, each hyper-elementary subgroup
of  ${\mathbb{Z}}_{q^s} \rtimes_\alpha {\mathbb{Z}}_{t_s}$ is conjugate to a subgroup H, such that one of the following is true:
\\
Case (1), $[{\mathbb{Z}}_{t_s},\pi(H)] \geq 4n^2$, where $\pi : {\mathbb{Z}}_{q^s} \rtimes_\alpha {\mathbb{Z}}_{t_s} \mapsto {\mathbb{Z}}_{t_s}$ is the natural epimorphism.\\
Case (2), H is a subgroup of $\{0\}\rtimes_\alpha{\mathbb{Z}}_{t_s}$ \\

We proceed under the assumption that case (1) holds; i.e. we assume $[ {\mathbb{Z}}_{t_s}, \pi(H)] = m \geq 4n^2$. Choose $E_H$ to be the real line $\mathbb{R}$ (dimension 1), with the simplicial structure with $m \mathbb{Z}$ as vertices. There is a $\Gamma$ action on $\mathbb{R}$ by $(x,k) y = y + k$ for $(x,k) \in {\mathbb{Z} {[{\frac{1}{d}} ]}} \rtimes_\alpha \mathbb{Z}$, $y \in \mathbb{R}$. Note $\Gamma$ does not act simplicially on $\mathbb{R}$, but $\bar{H}$ does. The stabilizers of a vertex or an edge is $\bar{H} \bigcap {\mathbb{Z} {[{\frac{1}{d}} ]}} \rtimes_\alpha \{0\} $, hence an abelian group. Define $f_H: \Gamma \rightarrow E_H$, by mapping $(x,k) \in {\mathbb{Z} {[{\frac{1}{d}} ]}} \rtimes_\alpha \mathbb{Z}$ to $ k \in \mathbb{R}$, $f_H$ is an $\bar{H}$ equivariant map. We need to prove for $g,h \in \Gamma$, if $d_\Gamma(g,h) < n$, then $d^1_{E_H} (f_H(g),f_H(h)) < \frac{1}{n}$.

Let $g = (x_1,k_1)$ and $h = (x_2, k_2)$, then $f_H(g) = k_1$ and $f_H(h) = k_2$. Hence $d(f_H(g),f_H(h)) = |k_1 - k_2|$. On the other hand, $d_\Gamma(g,h)= d(g(P_0,0), h(P_0,0)) = d(h^{-1}g(P_0,0),(P_0,0))$ (Here $P_0$ is the base point defined in section \ref{model}). Note $h^{-1}g = (d^{-k_2}(x_1 - x_2), k_1 - k_2)$, while $h^{-1}g$ acts on $(P_0, 0)$ as the matrix\\ \hspace*{1cm} $\left(
\begin{array}{cc}
    d^{-(k_1 - k_2)} & d^{k_2}(x_1 - x_2) \\
   0 & 1 \\
    \end{array}\right) $ = $\left(
\begin{array}{cc}
    1 & d^{k_2}(x_1 - x_2) \\
   0 & 1 \\
    \end{array}\right)$  $\left(
\begin{array}{cc}
    d^{-(k_1 - k_2)} & 0 \\
   0 & 1 \\
    \end{array}\right).$\\ By checking the action on $P_0$, one sees that $h^{-1}g$ moves $P_0$ at least distance $|k_1 - k_2|$ away. By Lemma \ref{dd},   $|k_1 - k_2| \leq d(h^{-1}g(P_0,0),(P_0,0)) = d(g,h)$. So $d(f_H(g),f_H(h)) <d_\Gamma(g,h)< n$. One can easily check now that $d^1_{E_H} (f_H(g),f_H(h)) < \frac{1}{n}$ since the simplicial structure we put on $\mathbb{R}$ has edge length greater than $m$ while $m > 4n^2$. Hence we have completed the proof in this case.\\

Now we proceed to Case(2)\label{marker}. In this case H is a subgroup of $\{0\}\rtimes_\alpha{\mathbb{Z}}_{t_s}$. Denote $\{0\}\rtimes_\alpha{\mathbb{Z}}_{t_s}$ by $K$, then  $\bar{K} = \Delta_n^{-1}(K) =  ({q^s ~\mathbb{Z} {[{\frac{1}{d}} ]}}) \rtimes_\alpha \mathbb{Z} $, $\bar{H} \subseteq \bar{K}$. Define a monomorphism $\varphi: \Gamma \rightarrow \Gamma$ by \\*
  $\hspace*{20mm} \varphi ((x,k)) =({q^s}x,k)$, for $(x,k) \in {\mathbb{Z} {[{\frac{1}{d}} ]}} \rtimes_\alpha \mathbb{Z} $. \\*
Note that its image is $\bar{K}$; hence ${\varphi}^{-1}:\bar{K} \rightarrow \Gamma $ is well defined. Note also that $\varphi$ can also be considered as the conjugation by $\left(
\begin{array}{cc}
    q^s & 0 \\
   0 & 1 \\
    \end{array}\right)$ since $\left(
\begin{array}{cc}
    q^s & 0 \\
   0 & 1 \\
    \end{array}\right) \left(
\begin{array}{cc}
    d^{-k} & x \\
   0 & 1 \\
    \end{array}\right) \\\left(
\begin{array}{cc}
    q^{-s} & 0 \\
   0 & 1 \\
    \end{array}\right) = \left(
\begin{array}{cc}
    d^{-k} & q^s x \\
   0 & 1 \\
    \end{array}\right) $. Define now a self homeomorphism $F_{q^s}: T_d  \times \mathbb{R} \rightarrow T_d \times \mathbb{R}$ by \\*
$\hspace*{20mm} F_{q^s} (z,w) = \left(
\begin{array}{cc}
    q^s & 0 \\
   0 & 1 \\
    \end{array}\right) (z,w)$, for $(z,w) \in T_d\times \mathbb{R}$\\*
where the action of $\left(
\begin{array}{cc}
    q^s & 0 \\
   0 & 1 \\
    \end{array}\right)$ on $(z,w) \in T_d\times \mathbb{R}$ is the diagonal action explained in section \ref{model}, Remark \ref{act}. In fact, $\left(
\begin{array}{cc}
    q^s & 0 \\
   0 & 1 \\
    \end{array}\right)$ acts on the tree $T_d$ as an isometry which fixes the line $L_0$ while it acts on $\mathbb{R}$ by $\left(
\begin{array}{cc}
    q^s & 0 \\
   0 & 1 \\
    \end{array}\right)~w = {q^s}w$ for $w \in \mathbb{R}$. It is easy to see that $F_{q^s}$ is $\varphi$-semi-equivariant since $(x,k) \in \Gamma$ acts on $T_d \times \mathbb{R}$ as the matrix $\left(
\begin{array}{cc}
    d^{-k} & x \\
   0 & 1 \\
    \end{array}\right)$ and $\left(
\begin{array}{cc}
    q^s & 0 \\
   0 & 1 \\
    \end{array}\right) \left(
\begin{array}{cc}
    d^{-k} & x \\
   0 & 1 \\
    \end{array}\right) =       \left(
\begin{array}{cc}
    d^{-k} & {q^s}x \\
   0 & 1 \\
    \end{array}\right)\left(
\begin{array}{cc}
    q^s & 0 \\
   0 & 1 \\
    \end{array}\right)$. Recall $F_{q^s}$ is $\varphi$-semi-equivariant if    $F_{q^s} (g(z,w))$ =
    $\varphi(g) F_{q^s}((z,w))$ for every $g \in \Gamma$.

Now choose $\lambda$ to be $n$, then by Proposition \ref{lte}, we have a $\mathcal{VC}yc$-cover $\mathcal{U}$ of dimension $N$ of $HFS(T_d \times\mathbb{R})$  satisfying the following. There exists $\delta >0$ which depends on $\lambda$ such that for any point $x \in HFS(T_d \times\mathbb{R})$, there exists $U_x \in \mathcal{U}$ such that $ (\Phi_{[-n, n]}(x))^{\delta}   ~\subseteq~  U_x$. We need to define an $\bar{H}$-equivariant map from $T_d \times \mathbb{R} $ to $HFS(T_d \times\mathbb{R})$. Using this map we can pull the $\mathcal{VC}yc$-cover $\mathcal{U}$ on $HFS(T_d \times\mathbb{R})$ back to $T_d \times \mathbb{R}$ and get a $\mathcal{VC}yc$-cover on $T_d \times \mathbb{R}$.  Consider the following diagram\\*

$\hspace*{10mm}\begin{CD}
\Gamma  @> {\eta} >> T_d \times \mathbb{R} @>{F_{q^s}^{-1}}>> T_d \times \mathbb{R} @> {\Psi_\tau}>>  HFS(T_d \times\mathbb{R}) \\
\end{CD}$ \\*
~\\*
where $\eta$ is the inclusion map given by the orbit of the base point, and $\Psi_\tau = \Phi_{\tau} \circ \Psi$ as defined in remark \ref{remfl}. In order to guarantee the composition is $\bar{H}$ equivariant, we change the action of $\bar{H}$ on the second $T_d \times \mathbb{R}$ (i.e., the image of $F_{q^s}^{-1}$)  by defining $ \bar{h} \bullet (z,w) = {\varphi}^{-1} (\bar{h})(z,w)$, for any $\bar{h} \in \bar{H}$, $(z,w) \in T_d\times \mathbb{R}$. We also define an $\bar{H}$ action on $HFS(T_d \times\mathbb{R})$ using the same method, by composing the action of $\Gamma$ with ${\varphi}^{-1}$.  Now using the composite $ {\Psi_\tau} \circ F_{q^s}^{-1}$, we pull back the cover on $HFS(T_d \times\mathbb{R})$ to $T_d \times \mathbb{R}$ and denote the nerve of this cover by $E(\lambda)$. It is a simplicial $\bar{H}$-complex of dimension $N$ whose stabilizers belong to $\mathcal{VC}yc$, which are either trivial or infinite cyclic. Also there is a canonical map from $T_d \times \mathbb{R}$ to $E(\lambda)$ (see for example, \cite {BLR2},section 4.1), denote this map as $\hat{f}_H$. Now for suitable choices of  $\tau$ and $q^s$, $f_H = \hat{f}_H \circ \eta$ will be the map used to make $\Gamma$ a Farrell-Hsiang group with respect to the family of abelian subgroups. In fact we choose $\tau = \ln {n} - \ln{\delta} +n$,  and rechoose $s$ so that $s>q$ and $ q^s >\bar{N}$, where $\bar{N}$ is determined by Lemma \ref{lpar} choosing $\epsilon = \frac{\delta^2}{2ne^n}$, $D = 2n$. Note that $d$ is fixed once our group $\Gamma$ is fixed.

Now we begin to prove that with these choices $f_H$ indeed works. Let $x_0,x_1 \subseteq T_d \times \mathbb{R}$ with $d(x_0,x_1) < n$, and let $x_2 =(q_1(x_0), q_2(x_1))$, recall that $q_1$, $q_2$ are projections from $T_d\times \mathbb{R}$ to $T_d$ or $\mathbb{R}$ respectively. By Lemma \ref{dd} and Corollary \ref{ddc}, we have $d(x_2,x_1) < d(x_0,x_1) < n$, $d(x_2,x_0) < 2d(x_0,x_1) < 2n$. Applying $F_{q^s}^{-1}$, by Lemma \ref{lpar} and our choice of $\epsilon$ and $D$, $d_{FS(X)}(\Psi_0(F_{q^s}^{-1}(x_0)),\Psi_0(F_{q^s}^{-1}(x_2))) < \epsilon = \frac{\delta^2}{2n{e^{n}}}$. On the other hand, $d(x_2,x_1) = d(F_{q^s}^{-1}(x_2),F_{q^s}^{-1}(x_1))$, which is the distance between $q_1(x_2)$ and $q_1(x_0)$ in the tree $T_d$ by Lemma \ref{dd}, since $x_2$ and $x_1$ has the same $\mathbb{R}$ coordinate. For convenience we will denote  $F_{q^s}^{-1}(x_i)$ by $\bar{x_i}$, $i = 0,1,2$, then  $d(\bar{x}_1,\bar{x}_2) < n$ and $d_{FS(X)}(\Psi_0(\bar{x}_0),\Psi_0(\bar{x}_2))) < \epsilon$.

\begin{figure}[h]

		\centering
		\includegraphics[width=1.0\textwidth]{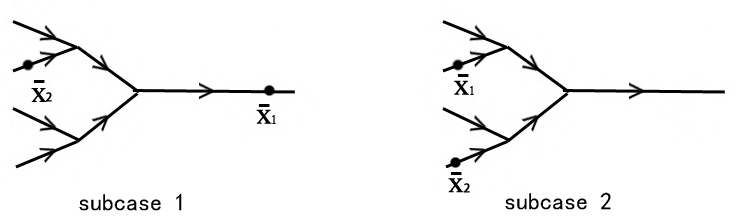}
		\centering
        \caption{Two subcases}\label{f2}
\end{figure}

As shown in Figure \ref{f2}  there are two subcases to consider depending on the positions of $\bar{x}_1 ~and~ \bar{x}_2$. Since $\tau = \ln{n} - \ln{\delta} +n$, by Lemma \ref{inq},
 $$d_{FS(X)}(\Psi_{\tau}({\bar{x}}_0), \Psi_{\tau}(\bar{x}_2))) \leq e^{|\tau|} d(\Psi_0(\bar{x}_0),\Psi_0(\bar{x}_2))) < \frac{{ne^{n}}}{\delta} \cdot \frac{\delta^2}{2n{e^{n}}} = \frac{\delta}{2}$$
In \textbf{subcase 1}, we first assume $\bar{x}_1$ is ahead of $\bar{x}_2$ as in Figure \ref{f2}. Let $\hat{d} = d(\bar{x}_1, \bar{x}_2)$, then there exists  an element of the cover $\mathcal{U}$, say $U$, such that $(\Psi_{[\tau - n ,\tau+n]}(\bar{x}_2))^\delta \subseteq U$ by Proposition \ref{lte}; hence  $\Psi_{\tau +\hat{d}}(\bar{x}_2) \in U$ since $\hat{d} <n$ .  On the other hand, as generalized geodesics, $\Psi_{\tau + \hat{d}}(\bar{x}_2)$ and $\Psi_{\tau}(\bar{x}_1)$ coincide for $t > -\tau$. In fact $\Psi_{\tau + \hat{d}}(\bar{x}_2)(t)$ maps $(-\infty, -\tau - \hat{d}]$ to the point $\bar{x}_2$ and $[-\tau -\hat{d}, \infty)$ isometrically to the geodesic starting with the point $\bar{x}_2$ and ending with the infinity point $\omega$; while $\Psi_{\tau}(\bar{x}_1)(t)$ maps $(-\infty, -\tau]$  to $\bar{x}_1$ and $[-\tau,\infty)$ isometrically to the geodesic starting with the point $\bar{x}_1$ and ending with the infinity point $\omega$. Hence

 $$d_{FS(X)}(\Psi_{\tau +\hat{d}}(\bar{x}_2),\Psi_\tau(\bar{x}_1)) = \int_{\mathbb{R}} \frac{d(\Psi_{\tau +\hat{d}}(\bar{x}_2)(t),\Psi_\tau(\bar{x}_1)(t))}{2e^{|t|}} dt \hspace*{44mm}$$

$$= \int_{(-\infty, -\tau-\hat{d}]} \frac{d(\Psi_{\tau +\hat{d}}(\bar{x}_2)(t),\Psi_\tau(\bar{x}_1)(t))}{2e^{|t|}} dt~~~~ ~+~~~~~\int_{[-\tau - \hat{d},-\tau]} \frac{d(\Psi_{\tau +\hat{d}}(\bar{x}_2)(t),\Psi_\tau(\bar{x}_1)(t))}{2e^{|t|}} dt \hspace*{20mm}$$

$$+  \int_{[-\tau, \infty]} \frac{d(\Psi_{\tau +\hat{d}}(\bar{x}_2)(t),\Psi_\tau(\bar{x}_1)(t))}{2e^{|t|}}dt~\hspace*{20mm}$$

$$< \int_{(-\infty, -\tau-\hat{d}]}\frac{d(\bar{x}_2, \bar{x}_1 )}{2e^{-t}}dt ~ +~   \int_{[-\tau-\hat{d}, -\tau]} \frac{d(\bar{x}_2, \bar{x}_1)}{2e^{-t}} dt        ~+ ~0~\hspace*{40mm}$$

$$=  \int_{(-\infty, -\tau-\hat{d}]}\frac{\hat{d}}{2e^{-t}}dt ~ +~   \int_{[-\tau-\hat{d}, -\tau]} \frac{\hat{d}}{2e^{-t}} dt \hspace*{62mm}$$
$$=\frac{1}{2}\hat{d} e^{-\tau}\leq ~ \frac{n}{2e^{\tau}} ~=~ \frac{~~n~~}{2e^{\ln{n} - \ln{\delta} +n}} = \frac{\delta}{2{e^{n}}} < \frac{\delta}{2}\hspace*{70mm}.$$


Therefore, we have proven that $\Psi_{\tau}(F_{q^s}^{-1}(x_0)),\Psi_{\tau}(F_{q^s}^{-1}(x_1)) \in (\Psi_{[\tau-n ,\tau+n]}(\bar{x}_2))^\delta \subseteq U$.  If $\bar{x}_2$ is ahead of $\bar{x}_1$, then the same calculation shows that
 $$d_{FS(X)}(\Psi_{\tau +\hat{d}}(\bar{x}_1),\Psi_\tau(\bar{x}_2)) < \frac{\delta}{2} $$
Consequently
 $$d_{FS(X)}(\Psi_{\tau +\hat{d}}(\bar{x}_1),\Psi_\tau(\bar{x}_0)) \leq  d_{FS(X)}(\Psi_{\tau +\hat{d}}(\bar{x}_1),\Psi_\tau(\bar{x}_2)) + d_{FS(X)}(\Psi_{\tau}({\bar{x}}_0), \Psi_{\tau}(\bar{x}_2))) $$
$$ < \frac{\delta}{2} + \frac{\delta}{2}= \delta \hspace*{15mm} $$
Therefore,$\Psi_{\tau}(F_{q^s}^{-1}(x_0)),\Psi_{\tau}(F_{q^s}^{-1}(x_1)) \in (\Psi_{[\tau-n ,\tau+n]}(\bar{x}_1))^\delta $ which will be contained in some member of the cover $\mathcal{U}$ by proposition \ref{lte}. Hence if we pull back the cover $\mathcal{U}$ via $\Psi_{\tau}\circ F_{q^s}^{-1}$, we get a cover $\mathcal{V}$ of $T_d \times \mathbb{R}$ and $x_0$ and $x_1$ will lie in the same member of $\mathcal{V}$.

 For \textbf{subcase 2}, note that $\Psi(\bar{x}_1)$ and $\Psi(\bar{x}_2)$ as generalized geodesics in $T_d \times \mathbb{R}$ will meet; suppose $y \in T_d \times \mathbb{R}$ is the first point where they meet. Let $d_i = d(\bar{x}_i,y)$ for $i=1,2$. Note $d_1 + d_2  < n$. We can assume $d_1 \geq d_2$ and let $\hat{d} = d_1 - d_2$, $\hat{d} < n$. Now we are almost in the
  same situation as \textbf{subcase 1}. The generalized geodesics, $\Psi_{\tau}(\bar{x}_1)$ and $\Psi_{\tau+ \hat{d}}(\bar{x}_2)$ will coincide for $t > -\tau + d_1$. And the same calculation shows that $d_{FS(X)}(\Psi_{\tau}(\bar{x}_1),\Psi_{\tau + \hat{d}}(\bar{x}_2)) < \frac{n}{2e^{\tau - d_1}} < \frac{n}{{2\frac{{ne^{n}}}{\delta}} e^{-d_1}} < \frac{\delta}{2}$. Therefore, \\
  $d_{FS(X)}(\Psi_{\tau +\hat{d}}(\bar{x}_0), \Psi_{\tau}(F_{q^s}^{-1}(x_1))) < \delta$, hence there exists an element $U$ of the cover $\mathcal{U}$ such that $\Psi_{\tau}(F_{q^s}^{-1}(x_0)),\Psi_{\tau}(F_{q^s}^{-1}(x_0)) \in (\Psi_{[\tau-n ,\tau+n]}(\bar{g}_0))^\delta \subseteq U$. Consequently if we pull back the cover $\mathcal{U}$, $x_0$ and $x_1$ will lie in the same member of this new cover $\mathcal{V}$ of $T_d \times  \mathbb{R}$.\\

So far, we have proved that  there exists a $\mathcal{VC}yc$-cover $\mathcal{V}$ of $T_d \times \mathbb{R}$, such that given any $x_0, x_1 \in T \times \mathbb{R}$ with $d(x_0,x_1) < n$, there exists a member of this cover containing both $x_0$ and $x_1$. In fact, our proof shows that any ball with radius less than $\frac{n}{2}$ will lie in some element of the cover. Hence we can apply the following lemma from \cite{BLR2}, proposition 5.3, page 47. \label{marker2}

\begin{lem}\label{cont}
Let $X = (X,d)$  be a metric space and $\beta \geq 1$. Suppose $\mathcal{U}$ is an open cover of $X$ of dimension less than or equal to $N$ with the following property:

$ For ~every~x\in X, there~ exists~U \in ~ \mathcal{U} ~ such ~ that ~ the ~ \beta ~ ball ~ around ~ x ~ lies ~ in ~ U.$
~\\*
Then if we denote the nerve of the cover as $N(\mathcal{U})$, the canonical map $\rho: X \rightarrow N(\mathcal{U}) $ has the following contracting property. If $d(x,y) \leq \frac{\beta}{4N}$, then
             $$d^1(\rho(x), \rho(y))   \leq  \frac{16N^2}{\beta} d(x,y).$$
\end{lem}

 By Proposition \ref{lte}, the dimension of the cover is less than a fixed number $N$ in our situation. For any $g_0,g_1 \in \Gamma$ such that $d(g_0,g_1) < n$, choose $\beta = 16N^2 n^2$, i.e., in the arguments of Case(2)  from page \pageref{marker} to page \pageref{marker2} before Lemma \ref{cont} , replace $n$ by $\bar{n} = 32N^2 n^2$, hence any ball of radius $16N^2 n^2$ will lie in the same element of the cover.  Then $d(g_0,g_1) < n < 4Nn^2= ~\frac{\beta}{4N}$, by Lemma \ref{cont}, $d^1(f_H(g_0), f_H(g_1))   \leq  \frac{16N^2}{\beta}~ d(g_0,g_1) = \frac{16N^2}{16 N^2 n^2} ~d(g_0,g_1) \leq \frac{1}{n^2} ~n = \frac{1}{n}$. Hence we have finished our proof.

\section{Further results} \label{fur}
In this section we extend our results on the Farrell-Jones conjecture to more general groups. As in Section \ref{gts}, we denote the unit group of a ring $R$ by $U(R)$.

Let $\hat{\Gamma}$ be the following matrix subgroup of $GL_2(\mathbb{Z}[\frac{1}{p}])$, where $p$ is a prime number \\*
$\hspace*{20mm} \hat{\Gamma} = \{~ \left(
\begin{array}{cc}
    \pm p^{-k} & x \\
   0 & 1 \\
    \end{array}\right) ~|~ k \in \mathbb{Z} ~and ~ x \in  \mathbb{Z}[\frac{1}{p}] ~\}$\\*
Note that the Baumslag-Solitar group $\Gamma = {\mathbb{Z} {[{\frac{1}{p}} ]}} \rtimes \mathbb{Z}$ we studied before, is a normal subgroup of index $2$ in $\hat{\Gamma}$
and that $\hat{\Gamma} = Aff(\mathbb{Z}[\frac{1}{p}])$; i.e. $\hat{\Gamma}$ is the full affine group of the ring $\mathbb{Z}[\frac{1}{p}]$. More precisely, $\hat{\Gamma}$ is isomorphic $ {\mathbb{Z} {[{\frac{1}{p}} ]}} \rtimes U(\mathbb{Z} {[{\frac{1}{p}} ]}) $ and $U(\mathbb{Z} {[{\frac{1}{p}} ]}) = \{ \pm p^{-k} ~|~ k \in \mathbb{Z}\} \cong \mathbb{Z} \oplus \mathbb{Z}_2$. Note that virtually cyclic groups are virtually abelian while the Farrell-Jones conjecture is true for virtually abelian groups. Hence by the transitivity principle, if we can prove $\hat{\Gamma}$ is a Farrell-Hsiang group with respect to the family of virtually abelian subgroups, then the Farrell-Jones conjecture is true for $\hat{\Gamma}$. By checking our proof for $\Gamma$, one sees that the only thing we need to take care of is section \ref{gts}.

For $q$ a prime number much bigger than $p$, define $U_{q^s}$ to be the subgroup of $U(\mathbb{Z}_{q^s})$ generated by $-1$ and $p$. Note since $q$ is an odd prime number, $U(\mathbb{Z}_{q^s})$ is a cyclic group (see for example \cite{KR}, chapter 4), hence $U_{q^s}$ is also a cyclic group, denote its generator by $d \in U(\mathbb{Z}_{q^s})$ and its order by $t'_s$. Then $t'_s = t_s$ or $2t_s$, where $t_s = {t_1} q^{k_s}$ is the order of $p$ in $U(\mathbb{Z}_{q^s})$ as defined in section \ref{gts}. Hence, $t'_1$, the order of $d$ in $\mathbb{Z}_{q}$, equals to either $t_1$ or $2t_1$, and $t'_s = {t'_1} q^{k_s}$. Note that $U_{q^s}$ as a subgroup of $U(\mathbb{Z}_{q^s})$ will have a canonical action on $\mathbb{Z}_{q^s}$ via multiplication, hence the semidirect product $\mathbb{Z}_{q^s} \rtimes U_{q^s}$ is well defined. And there is a quotient homomorphism $ \Delta_n :  {\mathbb{Z} {[{\frac{1}{p}} ]}} \rtimes U(\mathbb{Z} {[{\frac{1}{p}} ]}) \rightarrow \mathbb{Z}_{q^s} \rtimes U_{q^s}$.

\begin{lem}\label{gtr}
For any given integer $n >0$, let $q$  be a prime greater than $p^n$ and $s > q$, then every hyper-elementary subgroup H
of  $\mathbb{Z}_{q^s} \rtimes U_{q^s} = \mathbb{Z}_{q^s} \rtimes_\alpha {\mathbb{Z}}_{t'_s}$ is conjugate to a subgroup $\bar{H}$, such that one of the following is true:
\\
(1) the index $[\pi'(\bar{H}), {\mathbb{Z}}_{t'_s}] \geq n$, where $\pi' : {\mathbb{Z}}_{q^s} \rtimes_\alpha {\mathbb{Z}}_{t'_s} \mapsto {\mathbb{Z}}_{t'_s}$ is the natural epimorphism.\\
(2) $\bar{H}$ is a subgroup of $\{0\}\rtimes_\alpha{\mathbb{Z}}_{t'_s}$, and $q^s \geq n$.
\end{lem}
\begin{rem}
If we write $\mathbb{Z}_{q^s} \rtimes U_{q^s} $ as $ \mathbb{Z}_{q^s} \rtimes_\alpha {\mathbb{Z}}_{t'_s}$, then $\alpha$ is a multiplication by $d$, where $d \in U(\mathbb{Z}_{q^s})$ generates $U_{q^s}$.

\end{rem}

\Proof $\mathbb{Z}_{t_1}$ is a subgroup of index $1$ or $2$ in $\mathbb{Z}_{t'_1}$ depending on whether $t'_1 = t_1$ or $2 t_1$. Likewise $\mathbb{Z}_{t_s}$ is a subgroup of index $1$ or $2$ in $\mathbb{Z}_{t'_s}$ depending on whether $t'_1 = t_1$ or $2 t_1$.

Let $H' = H \cap \mathbb{Z}_{q^s} \rtimes {\mathbb{Z}}_{t_s}$, then $H'$ is a subgroup of index $1$ or $2$ of H. If the index is $1$, then we are done by Corollary \ref{c2}. Hence we can assume the index is $2$;  note that $H'$ is also a hyper-elementary group in this case.  Therefore by Corollary \ref{c2}, either

\begin{itemize}
\item [(1')] $[\pi'(H'), {\mathbb{Z}}_{t_s}] \geq n$ or
\item[(2')] $H'$ can be conjugate to a subgroup of  $\{0\}\rtimes_\alpha {\mathbb{Z}}_{t_s} \subseteq \mathbb{Z}_{q^s} \rtimes _\alpha {\mathbb{Z}}_{t_s}\subseteq \mathbb{Z}_{q^s} \rtimes_\alpha {\mathbb{Z}}_{t'_s}$.
\end{itemize}

It is clear that  (1') $ \Longrightarrow (1)$. Therefore we assume (2') occurs. And hence we can also (after this conjugation) assume that $H'$ is a subgroup of $\{0\}\rtimes_\alpha{\mathbb{Z}}_{t_s}$. Therefore the $Ker(\pi'|_H)$ is either $\{0\}$ or $\mathbb{Z}_2$. But the order of $Ker(\pi'|_H)$ has to divide $q$, hence $Ker(\pi'|_H) = \mathbb{Z}_2$ is impossible. Hence $\pi' :H \rightarrow \hat{H} := \pi'(H) \subseteq {\mathbb{Z}}_{t'_s} $ is an isomorphism. So H is a cyclic group with generator say $(a,b) \in \mathbb{Z}_{q^s} \rtimes_\alpha {\mathbb{Z}}_{t_s}$. Now note Lemma \ref{conj} implies if $t'_1 \nmid b$, then H can be conjugate to a subgroup of $\{0\}\rtimes_\alpha{\mathbb{Z}}_{t'_s}$. If $t'_1 ~|~ b$, then order of $b$ divides $q^r$. But the order of $b$ is the same as $|H|$, and $2~|~ |H|$, which is a contradiction.
\qed

With this Lemma, one can now check that our methods for proving the Farrell-Jones conjecture for $\Gamma$ can be extended to $\hat{\Gamma}$.

\begin{prop} \label{aff}
The K- and L-theoretic Farrell-Jones conjecture is true for $\hat{\Gamma} = Aff(\mathbb{Z} {[{\frac{1}{p}} ]})$.
\end{prop}


\begin{cor}
The K- and L-theoretic Farrell-Jones conjecture is true for every semi-direct product $\mathbb{Z} {[{\frac{1}{p}} ]} \rtimes C$, where $p$ is any prime number and C is any virtually cyclic group.
\end{cor}

In order to prove this corollary we need the following result from \cite{BFL}, Theorem 1.7 which generalizes \cite{FJ}, Proposition 2.2.
\begin{lem} \label{exp}
Let $1 \rightarrow K \rightarrow G \xrightarrow{\psi} Q \rightarrow 1$ be an exact sequence of groups. Suppose that the group Q and for any virtually cyclic subgroup $V \subseteq Q$ the group $\psi^{-1}(V)$ satisfies the K-theoretic Farrell-Jones conjecture. Then G satisfies the K-theoretic Farrell-Jones conjecture. The same is true for the L-theoretic Farrell-Jones conjecture.
\end{lem}~\\*
Proof of the Corollary. The group $C$ acts on $\mathbb{Z} {[{\frac{1}{p}} ]}$ by some representation (homomorphism)\\*
$\hspace*{25mm} \varphi : C \rightarrow Aut(\mathbb{Z} {[{\frac{1}{p}} ]}) = U(\mathbb{Z} {[{\frac{1}{p}} ]}) ~=~\{\pm p^{-k}~|~ k \in \mathbb{Z} \}.$\\*
Hence there is a canonical homomorphism \\*
$\hspace*{25mm} \bar{\varphi} : \mathbb{Z} {[{\frac{1}{p}} ]}  \rtimes_{\varphi} C \rightarrow \mathbb{Z} {[{\frac{1}{p}} ]}  \rtimes   U(\mathbb{Z} {[{\frac{1}{p}} ]})  $\\*
whose kernel $K$ is a subgroup of $C$ and therefore also virtually cyclic. Consider the following exact sequence\\*

$\hspace*{20mm}\begin{CD}
~~~~~~~~~~~~~~~~~1  @> >> K @>>> \mathbb{Z} {[{\frac{1}{p}} ]}  \rtimes_{\varphi} C @> \bar{\varphi} >>  I @>  >> 1 \\ \end{CD}$  \\*
~\\*
where $I = image (\bar{\varphi})$. To prove the corollary, we want to apply Lemma \ref{exp} to this exact sequence. First by Theorem 1.8 from \cite{BFL}, the group I satisfies the K- and L-theoretic Farrell-Jones conjecture since its overgroup $\mathbb{Z} {[{\frac{1}{p}} ]}  \rtimes   U(\mathbb{Z}{[{\frac{1}{p}} ]}) $ does by Proposition \ref{aff}. Therefore to complete the proof, by Lemma \ref{exp}, it remains to show that, for each virtually cyclic subgroup V of I, ${\bar{\varphi}}^{-1}(V)$ also satisfies the K- and L-theoretic Farrell-Jones conjecture. But this follows from Theorem 0.1 of \cite{BFL} since ${\bar{\varphi}}^{-1}(V)$ is a virtually poly-$\mathbb{Z}$ group. Note that ${\bar{\varphi}}^{-1}(V)$ is an extension of V by K, which are both virtually poly-$\mathbb{Z}$ groups, hence by Lemma 4.2 (v) in \cite{BFL}, ${\bar{\varphi}}^{-1}(V)$ is a virtually poly-$\mathbb{Z}$ group. \qed

\begin{rem}
The groups $\mathbb{Z} {[{\frac{1}{p}} ]} \rtimes C$ (where the image of $\varphi$ is infinite) forms an interesting subclass of nearly crystallographic groups as defined in [\cite{FL}, page 309]. The importance of the class of nearly crystallographic groups lies in Theorem 1.2 of \cite{FL}, where it was shown that the fibered isomorphism conjecture for the stable topological pseudo-isotopy functor is true for all virtually solvable groups provided it is true for the much smaller class consisting of all nearly crystallographic groups. And the truth of the fibered isomorphism conjecture for a torsion free group $G$ implies that $Wh(G) = 0$.
\end{rem}

\begin{center}
    {\bf APPENDIX}
  \end{center}
\appendix
\vspace{2mm}

Let $G_d$ be the group~ $\{{\left(
\begin{array}{cc}
    d^n \cdot \frac{s_1}{s_2} & b \\
   0 & 1 \\
    \end{array}\right)} ~|~ s_1,~s_2~ are~ nonzero~ integers~ coprime~ to ~d, ~b ~\in \mathbb{Q}\}$, $T_d$ be the oriented infinite, regular, $(d+1)$-valent tree with edge length $1$ as in Section \ref{model}. In this appendix we explain how the group $G_d$ acts on the tree $T_d$ where $d$ is a positive integer.

If $d$ is a prime, the action is well known as we explained in section \ref{model}. For more information, see for example \cite{Se}, Chapter II section 1.

We now assume $d$ is a power of a prime, hence $d = d_1^j$ for some prime $d_1$. Since $G_{d_1^j}$ is a subgroup of $G_{d_1}$, $G_{d_1^j}$ also acts on $T_{d_1}$. The following definition of Stallings' folding is taken from \cite{BMF}.

\begin{defn}
Let T be a G-tree (recall that this means, in particular, that there are no inversions).
Consider two edges $e_1$ and $e_2$ in T that are incident to a common vertex $v$. By
$\phi :e_1 \rightarrow e_2$ denote the linear homeomorphism fixing $v$. Then define an equivalence
relation " $\sim$ " on T as the smallest equivalence relation such that:\\
(i) $x \sim \phi(x)$, for all $x \in e_1$ , and\\
(ii) if $x \sim y$ and $g \in G$ then $g(x) \sim g(y)$.\\
The quotient space $T/\sim $ is a simplicial tree with a natural simplicial action of G. It
might happen that G acts with inversions on $T/\sim $, in which case we introduce
a new equivalence class of vertices to obtain a G-tree. Call the quotient map
$T \rightarrow T/\sim $ a fold.
\end{defn}

The key thing is that after the folding, G  has an induced action on the new tree.

Now we take the tree to be $T_p$ which is the oriented infinite regular (p+1)-valent tree we defined before; c.f., figure \ref{tree} in section \ref{model} for $T_2$. We first explain how does the action of  $G_{p^2}$ on $T_p$ induce an action on the tree $T_{p^2}$. Let $P_0P_1$ be an edge in the specified horizontal line $L_0$  of the tree $T_p$ such that $f_p(P_0) = 0$ and $f_p(P_1) = -1$; see section \ref{model} for the terminology here. The matrix ${\left(
\begin{array}{cc}
    1 & d^{-1} \\
   0 & 1 \\
    \end{array}\right)} $ fixes $P_{1}$ and acts on the $p$ edges that going towards $P_1$ by cyclic permutation. We will denote the image of $P_0$ under ${\left(
\begin{array}{cc}
    1 & d^{-1} \\
   0 & 1 \\
    \end{array}\right)} $ by $Q_1$. The folding map we are going to use is $\phi: P_0P_1 \rightarrow Q_1P_1$ which fixes $P_1$. Note that since ${\left(
\begin{array}{cc}
    1 & d^{-1} \\
   0 & 1 \\
    \end{array}\right)} $ acts as cyclic permutation on the $p$ edges going towards $P_1$, all the $p$ edges are identified after the folding. Therefore after the folding there will be $p^2$ edges going towards $P_0$. Moreover, $G_{p^2}$ is generated by matrices of the following three forms:
 \begin{itemize}
 \item
type I: ${\left(
\begin{array}{cc}
    p^{2n}  & 0 \\
   0 & 1 \\
    \end{array}\right)} $;
    \item
   type II: ${\left(
\begin{array}{cc}
    \frac{s_1}{s_2} & 0 \\
   0 & 1 \\
    \end{array}\right)} $, where $s_1,s_2$ are positive integers that are coprime to $p$;
    \item
    type III: ${\left(
\begin{array}{cc}
    1 & b \\
   0 & 1 \\
    \end{array}\right)} $, where $b \in \mathbb{Q}$
 \end{itemize}

\begin{figure}[h]
		\begin{center}
		\includegraphics[width=1.0\textwidth]{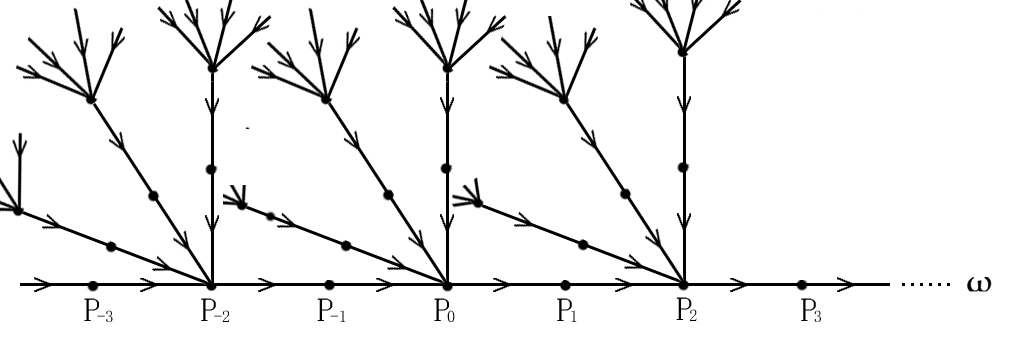}
		\end{center}
        \caption{$T_2$  after Stallings' folding}
        \label{rep}
	\end{figure}

Matrices of type I act as translation on $T_p$ which will change the Busmann function by an even number while matrices of type II and III leave the Busemann function unchanged. One sees now that the resulting tree is almost $T_{p^2}$ except it has some 2-valent vertices. We deleted these 2-valent vertices, and the group $G_{p^2}$ has an induced action on the new tree which is $T_{p^2}$. Figure \ref{rep} shows the resulting tree $T_2/\sim$, which is homeomorphic to $T_4$.
If we delete all the $2$-valent vertices in the resulting tree, it will be exactly $T_{4}$; in figure \ref{rep}, for example vertices $P_{2k+1}$ will be deleted.

One can further use this to get the action of $G_{p^l}$ on $T_{p^l}$ from the action of $G_{p^l}$ on $T_{p}$ by applying Stallings' folding for $l-1$ times.

If $d$ is not a power of a prime, let $d_1^{j_1} d_2^{j_2} \ldots d_m^{j_m}$ be its prime factorization. Note that $G_d$ is a subgroup of $G_{d_l^{j_l}}$ for $1 \leq l \leq m$; hence it acts on $T_{d_l^{j_l}}$, and therefore it will act on the product space $Y_d = T_{d_1^{j_1}} \times T_{d_2^{j_2}} \times \ldots \times T_{d_m^{j_m}}$ diagonally.  Now consider the ``diagonal subspace"
 $$\{y = (y_1,y_2,\ldots, y_{N})\in Y_{d}~ |~ f_{{d_l}^{j_l}}(y_l) = f_{{d_{l'}}^{j_{l'}}}(y_{l'}),~ for ~any ~1\leq l ~and~ l'\leq N \}$$
 where $f_{{d_l}^{j_l}}$ is the Busemann function on the corresponding $T_{d_{{d_l}^{j_l}}}$  we defined before; see Remark \ref{bus}. The diagonal space is an invariant subspace of $G_d$, hence $G_d$ has an induced action on it.  It is not too hard to show that this subspace is homeomorphic to the tree $T_d$ we defined before; we will think of them as the same. Therefore we have an induced action of $G_d$ on $T_d$. There is a natural Busemann function defined on $T_d$  by $f_{d}(y) = f_{d_1^{j_1}}(y_1)$ which is the same as we defined before; see Remark \ref{bus}.

~\\
Tom Farrell\\
DEPARTMENT OF MATHEMATICS, SUNY BINGHAMTON, NY,13902,U.S.A.\\
E-mail address: farrell@math.binghamton.edu\\
Xiaolei Wu\\
DEPARTMENT OF MATHEMATICS, SUNY BINGHAMTON, NY,13902,U.S.A.\\
E-mail address: xwu@math.binghamton.edu

\end{document}